\documentclass{gtart}

\usepackage{graphicx}
\usepackage{indentfirst}
\usepackage{amsmath,amsfonts,amsthm,amssymb}
\usepackage{mathrsfs}
\usepackage{amscd}

\newcommand{\homo}{{\mathrm H}}
\def\co{\colon\thinspace}
\DeclareMathAlphabet{\mathsfsl}{OT1}{cmss}{m}{sl}
\newcommand{\tensor}[1]{\mathsfsl{#1}}

\newtheorem{thm}{Theorem}[section]
\newtheorem{lem}[thm]{Lemma}
\newtheorem{cor}[thm]{Corollary}
\newtheorem{prop}[thm]{Proposition}
\newtheorem*{thm*}{Theorem}
\newtheorem*{claim}{Claim}

\theoremstyle{definition}
\newtheorem{defn}[thm]{Definition}
\newtheorem{rem}[thm]{Remark}
\newtheorem{construction}[thm]{Construction}

\begin{document}

\title{Knot Floer homology detects fibred knots}

\author{Yi Ni}
\address{Department of Mathematics, Princeton University\\ Princeton,  New Jersey 08544, USA
\\{\rm Current address:} Department of Mathematics, Columbia University\\Room 516, MC 440\/6, 2990 Broadway\\ New York, NY 10027, USA} \email{yni@math.columbia.edu}

\begin{abstract}
Ozsv\'ath and Szab\'o conjectured that knot Floer homology detects
fibred knots in $S^3$. We will prove this conjecture for
null-homologous knots in arbitrary closed 3--manifolds. Namely, if
$K$ is a knot in a closed 3--manifold $Y$, $Y-K$ is irreducible,
and $\widehat{HFK}(Y,K)$ is monic, then $K$ is fibred. The proof
relies on previous works due to Gabai, Ozsv\'ath--Szab\'o,
Ghiggini and the author. A corollary is that if a knot in $S^3$
admits a lens space surgery, then the knot is fibred.
\end{abstract}


\maketitle
\begin{center}
{\it Dedicated to Professor Boju Jiang on the occasion of his 70th
birthday}\end{center}

\section{Introduction}

Knot Floer homology was introduced independently by
Ozsv\'ath--Szab\'o \cite{OSz4} and by Rasmussen \cite{Ra}. For any
null-homologous knot $K\subset Y^3$ with Seifert surface $F$, one
can associate to it some abelian groups $\widehat{HFK}(Y,K,[F],i)$
for $i\in\mathbb Z$. The knot Floer homology
$$\widehat{HFK}(Y,K)\cong\oplus_{i\in\mathbb Z}\widehat{HFK}(Y,K,[F],i)$$
 is a finitely
generated abelian group.

A lot of topological information of the knot are contained in knot
Floer homology, in particular in the topmost filtration level. For
example, Ozsv\'ath and Szab\'o proved that the topmost filtration
level of $\widehat{HFK}$ for a knot in $S^3$ is exactly the genus
of the knot (see \cite{OSz5}).

When $K$ is a fibred knot, it is shown that the topmost group of
$\widehat{HFK}(Y,K)$ is a single $\mathbb Z$ \cite{OSz6}.
Ozsv\'ath and Szab\'o conjectured that the converse is also true
for knots in $S^3$ \cite{OSz9}.

In this paper, we are going to prove this conjecture. Our main
theorem is:
\begin{thm}\label{KnotFibre}
Suppose $K$ is a null-homologous knot in a closed, oriented,
connected 3--manifold $Y$, $Y-K$ is irreducible, and $F$ is a
genus $g$ Seifert surface of $K$. If
$\widehat{HFK}(Y,K,[F],g)\cong\mathbb Z$, then $K$ is fibred, and
$F$ is a fibre of the fibration.
\end{thm}

An oriented link $L$ in $Y$ is called a {\it fibred link}, if
$Y-L$ fibres over the circle, and $L$ is the oriented boundary of
the fibre. We have the following corollary of Theorem
\ref{KnotFibre}:

\begin{cor}\label{LinkFibre}
Suppose $Y$ is a closed, oriented, connected 3--manifold, $L$ is a
null-homologous oriented link in $Y$, $Y-L$ is irreducible, and
$F$ is a Seifert surface of $L$. If
$\widehat{HFK}\left(Y,L,\frac{|L|-\chi(F)}2\right)\cong\mathbb Z$,
then $L$ is a fibred link, and $F$ is a fibre of the fibration.
\end{cor}

The proof of this corollary will be given in Section 7.

A rational homology sphere $Y$ is called an {\it $L$--space}, if
the rank of $\widehat{HF}(Y)$ is equal to $|\homo_1(Y;\mathbb
Z)|$. Many 3--manifolds are $L$--spaces, for example, the
manifolds which admit spherical structures are $L$--spaces. An
immediate corollary of Theorem \ref{KnotFibre} is:

\begin{cor}
If a knot $K\subset S^3$ admits an $L$--space surgery, then $K$ is
a fibred knot. In particular, any knot that admits a lens space
surgery is fibred.
\end{cor}
\begin{proof}
As a corollary of \cite[Proposition 9.5]{OSz7}, if a rational
surgery on $K$ yields an $L$--space, then $K$ also admits an
$L$--space surgery with integer coefficient. Using
\cite[Theorem~1.2]{OSz5.5}, we conclude that
$\widehat{HFK}(K,g)\cong\mathbb Z$. Thus the desired result
follows from Theorem~\ref{KnotFibre}.
\end{proof}

\begin{cor}
Suppose $Y$ is an $L$--space, $K\subset Y$ is a null-homologous
knot with genus $g>1$. If the 0--surgery on $K$ is a surface
bundle over $S^1$, then $K$ is fibred.
\end{cor}
\begin{proof}
With the above conditions, one can prove that
$$HF^+(Y_0(K),[g-1])\cong\widehat{HFK}(Y,K,g).$$
In fact, the proof is exactly the same as the proof of
\cite[Corollary 4.5]{OSz4}, so we will not give the details. The
reader should note that, since $Y$ is an $L$--space, $HF^+(Y)$ is
isomorphic to the direct sum of some copies of $\mathbb
Z[U,U^{-1}]/U\mathbb Z[U]$, hence the map $\psi$ as in the proof
of \cite[Corollary 4.5]{OSz4} is surjective. Thus the argument
there can be used.

Now since $Y_0(K)$ fibres over the circle, we have
$HF^+(Y_0(K),[g-1])\cong\mathbb Z$, so
$\widehat{HFK}(Y,K,g)\cong\mathbb Z$. It is easy to show that $Y-K$
is irreducible, hence $K$ is fibred by Theorem~\ref{KnotFibre}.
\end{proof}

\begin{rem}
The homology class $[F]$ defines a homomorphism
$$f\co\pi_1(Y-K)\to\mathbb Z$$ by counting the intersection numbers
of $[F]$ with loops. The famous Stallings' Fibration Theorem
\cite{St} says that $K$ is a fibred knot with fibre in the
homology class of $[F]$ if and only if $\mathrm{ker}\:f$ is
finitely generated. Hence Theorem~\ref{KnotFibre} indicates a
mysterious relationship between Heegaard Floer homology and the
fundamental group.
\end{rem}

Theorem \ref{KnotFibre} was previously examined in various
sporadic cases, and some theoretical evidences were given in
\cite{Ni2}, but the first real progress was made by Ghiggini in
\cite{Gh}, where a strategy to approach this conjecture was
proposed, and the special case of genus-one knots in $S^3$ was
proved. Ghiggini's strategy plays an essential role in the present
paper, we will apply this strategy by using a method introduced by
Gabai \cite{G3}. Another key ingredient of this paper is the study
of decomposition formulas for knot Floer homology, which is based
on \cite{Ni2}.

The paper is organized as follows. In Section 2, we give some
backgrounds on sutured manifolds. We will introduce a sutured
manifold invariant which naturally comes from knot Floer homology.
We also present a construction of certain Heegaard diagrams. In
Section 3, we prove a homological version of the main theorem.
Section 4 is devoted to prove the horizontal decomposition formula
for knot Floer homology. In Section 5, we will prove a major
technical theorem: the decomposition formula for knot Floer
homology, in the case of decomposing along a separating product
annulus. In Section 6, we use Gabai's method to study Ghiggini's
strategy. As a result, we get a clearer picture of the sutured
manifold structure of the knot complement, namely,
Theorem~\ref{criterion}. Section 7 contains the proof of the main
theorem, we use the decomposition formulas we proved (especially
Theorem~\ref{SVerHFS}) to reduce the problem to the case that we
already know.

\noindent{\bf Acknowledgements.} This paper has been submitted to
Princeton University as the author's PhD thesis. We wish to thank
David Gabai and Zolt\'an Szab\'o for their guidance.

We would like to thank Paolo Ghiggini for many fruitful
discussions during the course of this work. This paper benefits a
lot from his work \cite{Gh}.

A version of Theorem \ref{criterion} was also proved by Ian Agol
via a different approach. We wish to thank him for some
interesting discussions.

We are grateful to Matthew Hedden, Andr\'as Juh\'asz, Tao Li,
Peter Ozsv\'ath, Jiajun Wang and Chenyang Xu for some helpful
conversations and their interests in this work. We are
particularly grateful to an anonymous referee for enormous
suggestions and corrections.

The author was partially supported by a Graduate School Centennial
Fellowship at Princeton University. Parts of the work were carried
out when the author visited UQAM and Peking University; he wishes
to thank Steve Boyer, Olivier Collin and Shicheng Wang for their
hospitality. The author extends his gratitude to the American
Institute of Mathematics and the Clay Mathematics Institute for
their subsequent support.

\section{Preliminaries}

\subsection{Sutured manifold decomposition}

The theory of sutured manifold decomposition was introduced by
Gabai in \cite{G1}. We will briefly review the basic definitions.

\begin{defn}
A {\it sutured manifold} $(M,\gamma)$ is a compact oriented
3--manifold $M$ together with a set $\gamma\subset \partial M$ of
pairwise disjoint annuli $A(\gamma)$ and tori $T(\gamma)$. The
core of each component of $A(\gamma)$ is a {\it suture}, and the
set of sutures is denoted by $s(\gamma)$.

Every component of $R(\gamma)=\partial M-\mathrm{int}(\gamma)$ is
oriented. Define $R_+(\gamma)$ (or $R_-(\gamma)$) to be the union
of those components of $R(\gamma)$ whose normal vectors point out
of (or into) $M$. The orientations on $R(\gamma)$ must be coherent
with respect to $s(\gamma)$, hence every component of $A(\gamma)$
lies between a component of $R_+(\gamma)$ and a component of
$R_-(\gamma)$.
\end{defn}

\begin{defn}\label{DefBal}\cite[Definition 2.2]{J}
A {\it balanced sutured manifold} is a sutured manifold
$(M,\gamma)$ satisfying\\ (1) $M$ has no closed components.\\
(2) Every component of $\partial M$ intersects $\gamma$
nontrivially.\\
(3)$\chi(R_+(\gamma))=\chi(R_-(\gamma))$.
\end{defn}

\begin{defn}
Let $S$ be a compact oriented surface with connected components
$S_1,\dots,S_n$. We define
$$x(S)=\sum_i\max\{0,-\chi(S_i)\}.$$
Let $M$ be a compact oriented 3--manifold, $A$ be a compact
codimension--0 submanifold of $\partial M$. Let $h\in\homo_2(M,A)$.
The {\it Thurston norm} $x(h)$ of $h$ is defined to be the minimal
value of $x(S)$, where $S$ runs over all the properly embedded
surfaces in $M$ with $\partial S\subset A$ and $[S]=h$.
\end{defn}

\begin{defn}
A sutured manifold $(M,\gamma)$ is {\it taut}, if $M$ is
irreducible, and $R(\gamma)$ is Thurston norm minimizing in
$\homo_2(M,\gamma)$.
\end{defn}

\begin{defn}
Let $(M,\gamma)$ be a sutured manifold, and $S$ a properly
embedded surface in M, such that no component of $\partial S$
bounds a disk in $R(\gamma)$ and no component of $S$ is a disk
with boundary in $R(\gamma)$. Suppose that for every component
$\lambda$ of $S\cap\gamma$, one of (1)--(3) holds:

(1) $\lambda$ is a properly embedded non-separating arc in
$\gamma$.

(2) $\lambda$ is a simple closed curve in an annular component $A$
of $\gamma$ in the same homology class as $A\cap s(\gamma)$.

(3) $\lambda$ is a homotopically nontrivial curve in a toral
component $T$ of $\gamma$, and if $\delta$ is another component of
$T\cap S$, then $\lambda$ and $\delta$ represent the same homology
class in $\homo_1(T)$.

Then $S$ is called a {\it decomposition surface}, and $S$ defines
a {\it sutured manifold decomposition}
$$(M,\gamma)\stackrel{S}{\rightsquigarrow}(M',\gamma'),$$
where $M'=M-\mathrm{int}(\mathrm{Nd}(S))$ and
\begin{eqnarray*}
\gamma'\;\;&=&(\gamma\cap M')\cup \mathrm{Nd}(S'_+\cap
R_-(\gamma))\cup
\mathrm{Nd}(S'_-\cap R_+(\gamma)),\\
R_+(\gamma')&=&((R_+(\gamma)\cap M')\cup S'_+)-\mathrm{int}(\gamma'),\\
R_-(\gamma')&=&((R_-(\gamma)\cap M')\cup
S'_-)-\mathrm{int}(\gamma'),
\end{eqnarray*}
where $S'_+$ ($S'_-$) is that component of
$\partial\mathrm{Nd}(S)\cap M'$ whose normal vector points out of
(into) $M'$.
\end{defn}

\begin{defn}
A decomposition surface is called a {\it product disk}, if it is a
disk which intersects $s(\gamma)$ in exactly two points. A
decomposition surface is called a {\it product annulus}, if it is
an annulus with one boundary component in $R_+(\gamma)$, and the
other boundary component in $R_-(\gamma)$.
\end{defn}

\begin{defn}
A decomposition surface $S$ in a balanced sutured manifold is
called a {\it horizontal surface}, if $S$ has no closed component,
$|\partial S|=|s(\gamma)|$,
$[S]=[R_+(\gamma)]\in\homo_2(M,\gamma)$, and $\chi(S)=\chi(R_+)$.
\end{defn}

\begin{defn}
A balanced sutured manifold $(M,\gamma)$ is {\it vertically
prime}, if any horizontal surface $S\subset M$ is parallel to
either $R_-(\gamma)$ or $R_+(\gamma)$.
\end{defn}

\subsection{Knot Floer homology and an invariant of sutured
manifolds}\label{SecHFS}

Heegaard Floer homology has been proved to have very close
relationship with Thurston norm \cite{OSz5}. The definition of
Thurston norm is purely topological (or combinatorial) \cite{T},
while the definition of Heegaard Floer homology involves
symplectic geometry and analysis. The bridge that connects these
two seemingly different topics is taut foliation.

The fundamental method of constructing taut foliations is sutured
manifold decomposition \cite{G1}. Thus we naturally expect that,
by studying the behavior of Heegaard Floer homology under sutured
manifold decomposition, we can get better understanding of the
relationship between Heegaard Floer homology and Thurston norm.

The first approach in such direction was taken in \cite{Ni2},
where the ``sutured Heegaard diagrams" for knots are introduced.
Using them, one can study a very special case of sutured manifold
decomposition: the Murasugi sum. Later on, Juh\'asz introduced an
invariant for sutured manifolds, called ``sutured Floer homology",
and proved some properties \cite{J}.

In this subsection, we will introduce another invariant of sutured
manifolds, which naturally comes from knot Floer homology. We will
prove some decomposition formulas for this invariant in Section 4
and Section 5.

Suppose $L\subset Y$ is a null-homologous oriented link, $F$ is a
Seifert surface of $L$. Decompose $Y-\mathrm{int}(\mathrm{Nd}(L))$
along $F$, we get a balanced sutured manifold $(M,\gamma)$. The
argument in \cite[Proposition 3.5]{Ni2} shows that, if we cut open
$Y$ along $F$, reglue by a homeomorphism of $F$ which is the
identity on the boundary, to get a new link $L'$ in a new manifold
$Y'$, then
$$\widehat{HFK}\left(Y,L,\frac{|\partial F|-\chi(F)}2\right)\cong
\widehat{HFK}\left(Y',L',\frac{|\partial F|-\chi(F)}2\right)$$ as
abelian groups. Therefore, $\widehat{HFK}\left(Y,L,\frac{|\partial
F|-\chi(F)}2\right)$ can be viewed as an invariant for the sutured
manifold $(M,\gamma)$. (For simplicity, let $\mathfrak i(F)$
denote $\frac{|\partial F|-\chi(F)}2$.)

More precisely, suppose $(M,\gamma)$ is a balanced sutured
manifold, $R_{\pm}(\gamma)$ are connected surfaces. There exists a
diffeomorphism
$$\psi\co R_+(\gamma)\to R_-(\gamma),$$
such that for each component $A$ of $\gamma$, $\psi$ maps one
boundary component of $A$ onto the other boundary component. We
glue $R_+(\gamma)$ to $R_-(\gamma)$ by $\psi$, thus get a manifold
with boundary consisting of tori. We fill each boundary torus by a
solid torus whose meridian intersects $s(\gamma)$ exactly once.
Now we get a closed 3--manifold $Y$. Let $L$ be the union of the
cores of the solid tori. The pair $(Y,L)$ is denoted by
$\iota(M,\gamma)$. Of course, $\iota(M,\gamma)$ depends on the way
we glue $R_+$ to $R_-$ and the way we fill in the solid tori. In
our case, changing the filling is equivalent to changing the
gluing map by Dehn twists along the components of $\partial R_+$.
By the remark in the last paragraph, the abelian group
$$\widehat{HFK}(\iota(M,\gamma),\mathfrak i(R_+(\gamma)))$$
is independent of the choice of the gluing, hence it is
independent of the choice of $\iota(M,\gamma)$.

\begin{prop}\label{WellHFS}
There is a well-defined invariant $\widehat{HFS}$ for balanced
sutured manifolds, which is characterized by the following two
properties.

(1) If $(M,\gamma)$ is a balanced sutured manifold,
$R_{\pm}(\gamma)$ are connected, then
$$\widehat{HFS}(M,\gamma)\cong\widehat{HFK}(\iota(M,\gamma),\mathfrak
i(R_+(\gamma))).$$

(2) If
$$(M,\gamma)\stackrel{a\times I}{\rightsquigarrow}(M',\gamma'),$$
is the decomposition along a product disk then
$$\widehat{HFS}(M,\gamma)\cong\widehat{HFS}(M',\gamma').$$
\end{prop}
\begin{proof}
The inverse operation of the decomposition along a product disk,
is ``adding a product 1--handle" with feet at the suture. We first
claim that, if $(M,\gamma)$ is a balanced sutured manifold with
$R_{\pm}(\gamma)$ connected, and $(M_1,\gamma_1)$ is obtained by
adding a product 1--handle to $(M,\gamma)$, then
$$\widehat{HFK}(\iota(M,\gamma),\mathfrak
i(R_+(\gamma)))\cong \widehat{HFK}(\iota(M_1,\gamma_1),\mathfrak
i(R_+(\gamma_1))).$$ In fact, one possible choice of
$\iota(M_1,\gamma_1)$ can be gotten by plumbing $\iota(M,\gamma)$
with $(S^2\times S^1,\Pi)$. Here $\Pi$ is a link in $S^2\times
S^1$ which consists of two copies of $\textsl{point}\times S^1$,
but with different orientations. Now we can apply
\cite[Lemma~4.4]{Ni2} to get the claim.

Given a balanced sutured manifold $(M,\gamma)$, we can add to it
some product 1--handles with feet at the suture to get a new
balanced sutured manifold $(M_1,\gamma_1)$, such that
$R_{\pm}(\gamma_1)$ are connected. We then define
$$\widehat{HFS}(M,\gamma)=\widehat{HFK}(\iota(M_1,\gamma_1),\mathfrak
i(R_+(\gamma_1))).$$

Now we want to prove that $\widehat{HFS}(M,\gamma)$ is independent
of the choice of $(M_1,\gamma_1)$. For this purpose, let
$(M_2,\gamma_2)$ be another sutured manifold obtained by adding
product 1--handles to $(M,\gamma)$, and $R_{\pm}(\gamma_2)$ are
connected. We can assume that the feet of the product 1--handles
for $(M_1,\gamma_1)$ and $(M_2,\gamma_2)$ are mutually different.
Let $(M_3,\gamma_3)$ be the sutured manifold obtained by adding
all these product 1--handles (either for $M_1$ or for $M_2$) to
$(M,\gamma)$. By the claim proved in the first paragraph, we have
\begin{eqnarray*}
\widehat{HFK}(\iota(M_3,\gamma_3),\mathfrak
i(R_+(\gamma_3)))&\cong&
\widehat{HFK}(\iota(M_1,\gamma_1),\mathfrak i(R_+(\gamma_1))),\\
\widehat{HFK}(\iota(M_3,\gamma_3),\mathfrak
i(R_+(\gamma_3)))&\cong&
\widehat{HFK}(\iota(M_2,\gamma_2),\mathfrak i(R_+(\gamma_2))).
\end{eqnarray*}
Therefore, $\widehat{HFS}(M,\gamma)$ is well-defined.

Property (1) holds by definition, and Property (2) can be proved
by the same argument as above.
\end{proof}

\subsection{Relative Morse functions and sutured diagrams}

Suppose $K$ is a null-homologous knot in $Y$, $F$ is a Seifert
surface for $K$. In \cite{Ni2}, the notion of ``sutured Heegaard
diagrams" was introduced. Such diagrams are useful to compute
$\widehat{HFK}(Y,K,[F],g)$. A construction of sutured Heegaard
diagrams was given in the proof of \cite[Theorem 2.1]{Ni2}.

In this subsection, we will present a slightly different
construction, which is based on relative Morse functions. This
construction will be useful later.

\begin{defn}\cite[Definition 2.2]{Ni2}
A double pointed Heegaard diagram
$$(\Sigma,\mbox{\boldmath${\alpha}$},
\mbox{\boldmath$\beta_0$}\cup\{\mu\},w,z)$$ for $(Y,K)$ is a {\it
sutured Heegaard diagram}, if it satisfies:

(Su0) There exists a subsurface $\mathcal P\subset\Sigma$, bounded
by two curves $\alpha_1\in\mbox{\boldmath${\alpha}$}$ and
$\lambda$. $g$ denotes the genus of $\mathcal P$.

(Su1) $\lambda$ is disjoint from \mbox{\boldmath$\beta_0$}. $\mu$
does not intersect any $\alpha$ curves except $\alpha_1$. $\mu$
intersects $\lambda$ transversely in exactly one point, and
intersects $\alpha_1$ transversely in exactly one point.
$w,z\in\lambda$ lie in a small neighborhood of $\lambda\cap\mu$,
and on different sides of $\mu$. (In practice, we often push $w,z$
off $\lambda$ into $\mathcal P$ or $\Sigma-\mathcal P$.)

(Su2) $(\mbox{\boldmath$\alpha$}-\{\alpha_1\})\cap\mathcal P$
consists of $2g$ arcs, which are linearly independent in
$\homo_1(\mathcal P,\partial\mathcal P)$. Moreover,
$\Sigma-\mbox{\boldmath$\alpha$}-\mathcal P$ is connected.
\end{defn}

\begin{construction}\label{RelMorse} Suppose $(M,\gamma)$ is the sutured
manifold obtained by cutting $Y-\mathrm{int}(\mathrm{Nd}(K))$ open
along $F$. Let $\psi\co R_+(\gamma)\to R_-(\gamma)$ be the gluing
map. Namely, if we glue $R_+(\gamma)$ to $R_-(\gamma)$ by $\psi$,
then we get back the manifold $Y-\mathrm{int}(\mathrm{Nd}(K))$. We
will construct a Heegaard diagram for the pair $(Y,K)$. The
construction consists of 4 steps.

\noindent{\bf Step 0.}\quad{\sl A relative Morse function}

Consider a self-indexed relative Morse function $u$ on
$(M,\gamma)$. Namely, $u$ satisfies:

(1) $u(M)=[0,3]$, $u^{-1}(0)=R_-(\gamma)$,
$u^{-1}(3)=R_+(\gamma)$.

(2) $u$ has no degenerate critical points. $u$ is the standard
height function near $\gamma$. $u^{-1}\{\textrm{critical points of
index $i$}\}=i$.

(3) $u$ has no critical points on $R(\gamma)$.

Let $\widetilde F=u^{-1}(\frac32)$. $\partial\widetilde F$ is
denoted by $\widetilde{\lambda}$. Similarly, the boundary
components of $R_{\pm}(\gamma)$ are denoted by $\lambda_{\pm}$.

Suppose $u$ has $r$ index--$1$ critical points, then the genus of
$\widetilde F$ is $g+r$. The gradient $-\nabla u$ generates a flow
$\phi_t$ on $M$. There are $2r$ points on $R_+(\gamma)$, which are
connected to index--$2$ critical points by flowlines. We call
these points ``bad" points. Similarly, there are $2r$ bad points
on $R_-(\gamma)$, which are connected to index--$1$ critical
points by flowlines.

\noindent{\bf Step 1.}\quad{\sl Construct the curves}

Choose a small disk $D_+$ in a neighborhood of $\lambda_+$ in
$R_+(\gamma)$. Choose an arc $\delta_+\subset R_+(\gamma)$
connecting $D_+$ to $\lambda_+$. Flow $D_+$ and $\delta_+$ by
$\phi_t$, their images on $\widetilde F$ and $R_-(\gamma)$ are
$\widetilde D,D_-,\widetilde{\delta},\delta_-$. (Of course, we
choose $D_+$ and $\delta_+$ generically, so that the flowlines
starting from them do not terminate at critical points.) We can
suppose the gluing map $\psi$ maps $\delta_+$ onto $\delta_-$,
$D_+$ onto $D_-$. Let
$A_{\pm}=R_{\pm}(\gamma)-\mathrm{int}(D_{\pm})$, $\widetilde
A=\widetilde{F}-\mathrm{int}(\widetilde{D})$.

On $\widetilde F$, there are $r$ simple closed curves
$\widetilde{\alpha}_{2g+2},\dots,\widetilde{\alpha}_{2g+1+r}$,
which are connected to index--$1$ critical points by flowlines.
And there are $r$ simple closed curves
$\widetilde{\beta}_{2g+2},\dots,\widetilde{\beta}_{2g+1+r}$, which
are connected to index--$2$ critical points by flowlines.

Choose $2g$ disjoint arcs $\xi^-_2,\dots,\xi^-_{2g+1}\subset A_-$,
such that their endpoints lie on $\lambda_-$, and they are
linearly independent in $\homo_1(A_-,\partial A_-)$. We also
suppose they are disjoint from $\delta_-$ and the bad points. Let
$\xi^+_i=\psi^{-1}(\xi^-_i)$. We also flow back
$\xi^-_2,\dots,\xi^-_{2g+1}$ by $\phi_{-t}$ to $\widetilde F$, the
images are denoted by
$\widetilde{\xi}_2,\dots,\widetilde{\xi}_{2g+1}$.

Choose $2g$ disjoint arcs $\eta^+_2,\dots,\eta^+_{2g+1}\subset
A_+$, such that their endpoints lie on $\partial D_+$, and they
are linearly independent in $\homo_1(A_+,\partial A_+)$. We also
suppose they are disjoint from $\delta_+$ and the bad points. Flow
them by $\phi_t$ to $\widetilde F$, the images are denoted by
$\widetilde{\eta}_2,\dots,\widetilde{\eta}_{2g+1}$.

\noindent{\bf Step 2.}\quad{\sl Construct a diagram}

Let
$$\Sigma=A_+\cup \widetilde A\cup\{\textrm{2 annuli}\},$$
where one annulus is $\phi_{[0,3/2]}(\partial D_+)$, the other is
$\phi_{[0,3/2]}(\widetilde{\lambda})$. Let
$$\alpha_i=\xi^+_i\cup\widetilde{\xi}_i\cup\{\textrm{2 arcs}\},$$
where the 2 arcs are vertical arcs connecting $\xi^+_i$ to
$\widetilde{\xi}_i$ on an annulus, $i=2,\dots,2g+1$. Similarly,
let
$$\beta_i=\eta^+_i\cup\widetilde{\eta}_i\cup\{\textrm{2 arcs}\}.$$
Let $\alpha_1=\partial\widetilde D$,
$\lambda=\widetilde{\lambda}$,
$$\mu=\delta_+\cup\widetilde{\delta}\cup\{\textrm{2 arcs}\}.$$

Let
\begin{eqnarray*}
\mbox{\boldmath${\alpha}$} &= &
\{\alpha_2,\dots,\alpha_{2g+1}\}\cup\{\widetilde{\alpha}_{2g+2},\dots,\widetilde{\alpha}_{2g+1+r}\}\cup\{\alpha_1\},\\
\mbox{\boldmath${\beta}$} &= &
\{\beta_2,\dots,\beta_{2g+1}\}\cup\{\widetilde{\beta}_{2g+2},\dots,\widetilde{\beta}_{2g+1+r}\}\cup\{\mu\}.
\end{eqnarray*}
Pick two base points $w,z$ near $\lambda\cap\mu$, but on different
sides of $\mu$.

\noindent{\bf Step 3.}\quad{\sl Check that
$(\Sigma,\mbox{\boldmath${\alpha}$}, \mbox{\boldmath$\beta$},w,z)$
is a Heegaard diagram for $(Y,K)$}

This step is quite routine, we leave the reader to check the
following

\noindent(A) $\Sigma$ separates $Y$ into two genus--$(2g+1+r)$
handlebodies $U_1,U_2$, where
$$U_1=u^{-1}[0,\frac32]\bigcup\phi_{[0,\frac32]}(D_+),\quad U_2=u^{-1}[\frac32,3]-\phi_{[0,\frac32]}(D_+).$$
Every curve in
$\mbox{\boldmath${\alpha}$}$ bounds a disk in $U_1$, every curve
in $\mbox{\boldmath${\beta}$}$ bounds a disk in $U_2$.

\noindent(B) $\Sigma-\mbox{\boldmath${\alpha}$}$ is connected,
$\Sigma-\mbox{\boldmath${\beta}$}$ is connected.

\noindent(C) The two base points $w,z$ give the knot $K$ in $Y$.

Then
$$(\Sigma,\mbox{\boldmath${\alpha}$},
\mbox{\boldmath$\beta$},w,z)$$ is a Heegaard diagram for $(Y,K)$.
It is not hard to see that this is a sutured Heegaard diagram.
\qed
\end{construction}

The idea of considering a relative Morse function on the sutured
manifold appeared in \cite{Go}. This idea was communicated to the
author by Zolt\'an Szab\'o, and was used by Andr\'as Juh\'asz to
define a sutured manifold invariant \cite{J}.

\begin{lem}\label{admissible}
We can wind the $\widetilde{\alpha}$--curves in $\widetilde
A-\widetilde{\delta}$, and wind the $\xi^+$--curves in
$A_+-\delta_+$, to get a weakly admissible Heegaard diagram.
\end{lem}
\begin{proof}
We claim that,
$$E=\widetilde{A}-\cup_{i=2}^{2g+1}\widetilde{\xi}_i-\cup_{j=2g+2}^{2g+1+r}\widetilde{\alpha}_j-\widetilde{\delta}$$
is connected. In fact, since $A_-$ is connected, and $\xi^-_j$'s
are linearly independent in $\homo_1(A_-,\partial A_-)$, it is
easy to see that $E_-=A_--\cup_{i=2}^{2g+1}{\xi}^-_i-{\delta}_-$
is connected. $E$ is obtained by removing $2r$ disks from $E_-$,
so $E$ is also connected.

Since $E$ is connected, we can find simple closed curves
$$\theta_{2g+2},\dots,\theta_{2g+1+r}\subset\widetilde A-\cup\widetilde{\xi}_i-\widetilde{\delta},$$ which
are {\it geometrically dual} to
$\widetilde{\alpha}_{2g+2},\dots,\widetilde{\alpha}_{2g+1+r}$.
Namely, $\theta_j$ is disjoint from $\widetilde{\alpha}_i$ when
$j\ne i$, and $\theta_i$ intersects $\widetilde{\alpha}_i$
transversely at exactly one point.

We can also find closed curves in $A_+-\delta_+$, which are
geometrically dual to the $\xi^+$--curves. Now our desired result
follows from the argument in \cite[Proposition 3.3]{Ni2}.
\end{proof}

\section{A homological version of the main theorem}

\begin{prop}\label{HomoProd}
Suppose $K\subset Y$ is a null-homologous knot with Seifert surface
$F$, $(M,\gamma)$ is the sutured manifold obtained by cutting
$Y-\mathrm{int}(\mathrm{Nd}(K))$ open along $F$.

If $\widehat{HFK}(Y,K,[F],g)\cong\mathbb Z$, then $M$ is a
homology product, namely,
$$\homo_*(M,R_-(\gamma);\mathbb Z)\cong\homo_*(M,R_+(\gamma);\mathbb Z)\cong0.$$
\end{prop}

Fix a field $\mathbb F$, let $\mathit 0\in\mathbb F$ be the zero
element.

\begin{lem}\label{H_2(M)}
If $\widehat{HFK}(Y,K,[F],g)\cong\mathbb Z$, then
$\homo_2(M;\mathbb F)=0$.
\end{lem}
\begin{proof}
We use the sutured diagram constructed in Construction
\ref{RelMorse}. For the generators of $\widehat{HFK}(Y,K,[F],-g)$,
the intersection point $\mu\cap\alpha_1$ is always chosen. The
generators are supported in $\widetilde{A}$ \cite[Theorem
5.1]{OSz4}, hence the determinant of the $(2g+r)\times(2g+r)$
matrix
$$\tensor{V}=\begin{pmatrix}
\widetilde{\xi}_i\cdot\widetilde{\eta}_j
&\widetilde{\xi}_i\cdot\widetilde{\beta}_l\\
\widetilde{\alpha}_k\cdot\widetilde{\eta}_j
&\widetilde{\alpha}_k\cdot\widetilde{\beta}_l\\
\end{pmatrix}\quad{\begin{array}{l}2\le i,j\le2g+1\\2g+2\le k,l\le2g+1+r\end{array}}
$$
computes the Euler characteristic of
$\widehat{HFK}(Y,K,[F],g)\cong\mathbb Z$. Hence
$$\det\tensor{V}=\pm1.$$

Now if $\homo_2(M;\mathbb F)\ne0$, then there is a nontrivial
$\mathbb F$--linear combination
$$\gamma=\sum a_k\widetilde{\alpha}_k+\sum b_l\widetilde{\beta}_l,$$
which is 0 in $\homo_1(\widetilde{F};\mathbb F)$. So
$\gamma\cdot\widetilde{\alpha}_k=\mathit0$,
$\gamma\cdot\widetilde{\xi}_i=\mathit0$.

Since $\widetilde{\alpha}$--curves are linearly independent in
$\homo_1(\widetilde{\Sigma};\mathbb F)$, the coefficients $b_l$'s
are not all $\mathit0$. Without loss of generality, we can assume
$b_{2g+2}\ne\mathit0$.

Since the $\widetilde{\alpha}$--curves and
$\widetilde{\xi}$--curves are mutually disjoint, we have
$$\widetilde{\alpha}_k\cdot\sum b_l\widetilde{\beta}_l=\mathit0,\quad\widetilde{\xi}_i\cdot\sum
b_l\widetilde{\beta}_l=\mathit0.$$ So by elementary column
operations, we can change $\tensor{V}$ into a matrix with a zero
column, while the determinant of this new matrix is
$b_{2g+2}\det\tensor{V}\ne\mathit0$. This gives a contradiction.
\end{proof}

\begin{lem}\label{InclInj}
If $\widehat{HFK}(Y,K,[F],g)\cong\mathbb Z$, then the map
$$i_*:\homo_1(R_-(\gamma),\partial R_-(\gamma);\mathbb F)\to\homo_1(M,\gamma;\mathbb F)$$
is injective.
\end{lem}
\begin{proof}
$\homo_1(R_-(\gamma),\partial R_-(\gamma))$ is generated by the
$\xi^-$--curves. If $i_*$ is not injective, then there exists a
nontrivial linear combination
$$\zeta=\sum c_i\widetilde{\xi}_i,$$
which is homologous to a linear combination
$$\gamma=\sum a_k\widetilde{\alpha}_k+\sum b_l\widetilde{\beta}_l$$
in $\homo_1(\widetilde{F},\partial\widetilde{F};\mathbb F)$. We
have $(\zeta-\gamma)\cdot\widetilde{\beta}_l=\mathit0$. The
$\widetilde{\eta}$--curves are non-proper arcs in $\widetilde{F}$,
but we can connect the two endpoints of $\widetilde{\eta}_j$ by an
arc in $\widetilde{D}$ to get a closed curve in $\widetilde{F}$.
$\widetilde{\xi}_i,\widetilde{\alpha}_k,\widetilde{\beta}_l$ do
not intersect $\widetilde D$, so
$(\zeta-\gamma)\cdot\widetilde{\eta}_j=\mathit0$. Since
$\widetilde{\beta}$--curves and $\widetilde{\eta}$--curves are
mutually disjoint, we have
$$\widetilde{\eta}_j\cdot(\zeta-\sum a_k\widetilde{\alpha}_k)=\mathit0,\quad\widetilde{\beta}_l\cdot(\zeta-\sum a_k\widetilde{\alpha}_k)=\mathit0.$$
Now we can get a contradiction as in the proof of Lemma
\ref{H_2(M)}.
\end{proof}

\begin{proof}[Proof of Proposition \ref{HomoProd}]
 Let $R_{\pm}=R_{\pm}(\gamma)$. By
Lemma~\ref{InclInj}, we have the exact sequence
$$0\to\homo_2(R_-,\partial R_-;\mathbb F)\to\homo_2(M,\partial R_-;\mathbb F)\to\homo_2(M,R_-;\mathbb F)\to0.$$
Compare this exact sequence with
$$\homo_2(M;\mathbb F)\to\homo_2(M,\gamma;\mathbb F)\to\homo_1(\gamma;\mathbb F)\to\cdots.$$
Note that $$\homo_2(M,\gamma)=\homo_2(M,\partial
R_-),\quad\homo_2(R_-,\partial
R_-)\cong\homo_1(\gamma)\cong\mathbb Z,$$ and $\homo_2(M;\mathbb
F)=0$ by Lemma \ref{H_2(M)}, so
$$\homo_2(M,R_-;\mathbb F)=0.$$

By Poincar\'e duality, we have $$\homo_1(M,\partial M;\mathbb
F)\cong\homo^2(M;\mathbb F)\cong\homo_2(M;\mathbb F)\cong0.$$ So
we have the exact sequence
$$0\to\homo_2(M,\partial M;\mathbb F)\to\homo_1(\partial M,R_-;\mathbb F)\to\homo_1(M,R_-;\mathbb F)\to0.$$
Compare this exact sequence with
$$0\to\homo_1(R_-;\mathbb F)\to\homo_1(M;\mathbb F)\to\homo_1(M,R_-;\mathbb F)\to0.$$
Note that
$$\homo_1(\partial M,R_-;\mathbb F)\cong\homo_1(R_-;\mathbb F)\cong\mathbb
F^{2g},$$ $$\homo_1(M;\mathbb F)\cong\homo^1(M;\mathbb
F)\cong\homo_2(M,\partial M;\mathbb F),$$ we should have
$$\homo_1(M,R_-;\mathbb F)=0.$$

Hence we have proved that $\homo_*(M,R_-;\mathbb F)=0$ for any
field $\mathbb F$. So $$\homo_*(M,R_-;\mathbb Z)=0.$$ Similarly,
we have $\homo_*(M,R_+;\mathbb Z)=0$.
\end{proof}

\section{Horizontal decomposition}

\begin{thm}\label{HorKnot}
Let $K'\subset Y'$, $K''\subset Y''$ be two null-homologous knots.
Suppose $F',F''$ are two genus--$g$ Seifert surfaces for $K',K''$,
respectively. We construct a new manifold $Y$ and a knot $K\subset
Y$ as follows. Cut open $Y',Y''$ along $F',F''$, we get sutured
manifolds $(M',\gamma')$, $(M'',\gamma'')$. Now glue
$R_+(\gamma')$ to $R_-(\gamma'')$, glue $R_+(\gamma'')$ to
$R_-(\gamma')$, by two diffeomorphisms. We get a manifold $Z$ with
torus boundary. There is a simple closed curve $\mu\subset\partial
Z$, which is the union of the two cut-open meridians of $K',K''$.
We do Dehn filling along $\mu$ to get the manifold $Y$, the knot
$K$ is the core of the filled-in solid torus.

Our conclusion is
$$\widehat{HFK}(Y,K,[F'],g)\cong\widehat{HFK}(Y',K',[F'],g)\otimes\widehat{HFK}(Y'',K'',[F''],g),$$
as linear spaces over any field $\mathbb F$.
\end{thm}

\begin{rem}
We did not specify the gluings, since they will not affect our
result, thanks to \cite[Proposition 3.5]{Ni2}.
\end{rem}

\begin{rem}
We clarify some convention we are going to use throughout this
paper. A holomorphic disk in the symmetric product is seen as an
immersed subsurface of the Heegaard surface $\Sigma$. Suppose $Q$
is a subsurface of $\Sigma$, $D_1,\dots,D_n$ are the closures of
the components of $Q-\cup\alpha_i-\cup\beta_j$, choose a point
$z_k$ in the interior of $D_k$ for each $k$. If $\Phi$ is a
holomorphic disk, then $\Phi\cap Q$ denotes the immersed surface
$\sum_kn_{z_k}(\Phi)D_k$.
\end{rem}

\begin{proof}[Proof of Theorem \ref{HorKnot}]
The proof uses the techniques from \cite{Ni2}. We construct a
sutured Heegaard diagram
$(\Sigma',\mbox{\boldmath${\alpha'}$},\mbox{\boldmath${\beta'}$},w',z')$
for $(Y',K')$, as in the proof of \cite[Theorem 2.1]{Ni2}. The
reader may refer to Figure~1 there for a partial picture.

As a result, $\Sigma'$ is the union of two compact surfaces
$A',B'$, where $A'$ is a genus $g$ surface with two boundary
components $\alpha_{1A}',\lambda_A'$, and $B'$ is a genus $g+r'$
surface with two boundary components $\alpha_{1B}',\lambda_B'$.
$A'$ and $B'$ are glued together, so that $\alpha_{1A}'$ and
$\alpha_{1B}'$ become one curve $\alpha_1'$,
 $\lambda_A'$ and $\lambda_B'$ become one curve $\lambda'$.

We have
\begin{eqnarray*}
\mbox{\boldmath${\alpha'}$}&=&\{\alpha_1',\alpha_2',\dots,\alpha_{2g+1}',{\alpha}_{2g+2}',\dots,{\alpha}_{2g+1+r'}'\},\\
\mbox{\boldmath${\beta'}$}&=&\{\mu',\beta_2',\beta_3',\dots,\beta_{2g+1+r'}'\}.
\end{eqnarray*}
Here $\alpha_i'$ is the union of two arcs $\xi_i'\subset A'
,\overline{\xi'_i}\subset B'$, for $i=2,\dots,2g+1$. ${\alpha}_j'$
lie in $B'$, for $j=2g+2,\dots,2g+1+r'$. $\mu'$ is the union of
two arcs $\delta'\subset A'$, $\overline{\delta'}\subset B'$.
$\mu'$ intersects $\alpha_1'$ transversely in one point, and is
disjoint from all other $\alpha'$--curves. $\beta'_i$'s are
disjoint from $\lambda'$.

Similarly, we construct a sutured diagram
$(\Sigma'',\mbox{\boldmath${\alpha''}$},\mbox{\boldmath${\beta''}$},w'',z'')$.
$\Sigma''$ is the union of $A'',B''$. And the corresponding curves
are denoted by $\alpha_i'',\beta_i'',\dots$

Now we glue $A',B',A'',B''$ together, so that $\alpha_{1A}'$ and
$\alpha_{1B}'$ become one curve $\gamma_1'$, $\lambda_B'$ and
$\lambda_A''$ become one curve $\lambda'$, $\alpha_{1A}''$ and
$\alpha_{1B}''$ become one curve $\gamma_1''$, $\lambda_B''$ and
$\lambda_A'$ become one curve $\lambda''$. $\overline{\xi_i'}$ and
$\xi_i''$ are glued together to be a curve $\gamma_i'$,
$\overline{\xi_i''}$ and $\xi_i'$ are glued together to be a curve
$\gamma_i''$, $i=2,\dots,2g+1$. $\gamma_j'=\alpha_j'$ when
$j=2g+2,\dots,2g+1+r'$, $\gamma_k''=\alpha_k''$ when
$k=2g+2,\dots,2g+1+r''$. $\beta_i'$ and $\beta_i''$ are as before.
$\delta',\overline{\delta'},\delta'',\overline{\delta''}$ are
glued together to a closed curve $\mu$. We also pick two
basepoints $w,z$ near $\lambda'\cap\mu$, but on different sides of
$\mu$.

Let
\begin{eqnarray*}
\Sigma&=&A'\cup B'\cup A''\cup B''\\
\mbox{\boldmath${\gamma}$}&=&\{\gamma_1',\gamma_2',\dots,\gamma_{2g+1+r'}',\gamma_2'',\dots,\gamma_{2g+1+r''}''\},\\
\mbox{\boldmath${\beta}$}&=&\{\mu,\beta_2',\dots,\beta_{2g+1+r'}',\beta_2'',\dots,\beta_{2g+1+r''}''\}.
\end{eqnarray*}
Then
$(\Sigma,\mbox{\boldmath${\gamma}$},\mbox{\boldmath${\beta}$},w,z)$
is a Heegaard diagram for $(Y,K)$.

As in the proof of  \cite[Proposition 3.3]{Ni2}, we can wind
$\xi'_2,\dots,\xi'_{2g+1}$ in $A'-\delta'$,
$\gamma_{2g+2}',\dots,\gamma_{2g+1+r'}'$ in
$B'-\overline{\delta'}$, $\xi''_2,\dots,\xi''_{2g+1}$ in
$A''-\delta''$, $\gamma_{2g+2}'',\dots,\gamma''_{2g+1+r''}$ in
$B''-\overline{\delta''}$, so that the diagrams
$$(\Sigma,\mbox{\boldmath${\gamma}$},\mbox{\boldmath${\beta}$},w,z),(\Sigma',\mbox{\boldmath${\alpha'}$},\mbox{\boldmath${\beta'}$},w',z'),
(\Sigma'',\mbox{\boldmath${\alpha''}$},\mbox{\boldmath${\beta''}$},w'',z'')$$
become admissible, and any nonnegative relative periodic domain in
$A'$ or $A''$ (for these diagrams) is supported away from
$\lambda''$, $\lambda'$.

\begin{claim} If $\mathbf x$ is a generator of $\widehat{CFK}(Y,K,-g)$,
then $\mathbf x$ is supported outside $\mathrm{int}(A')\cup
\mathrm{int}(A'')$.
\end{claim}

Let $Y_0(K)$ be the manifold obtained from $Y$ by $0$--surgery on
$K$, $\underline{\mathfrak s}_{w,z}(\mathbf y)\in
\textrm{Spin}^c(Y_0)$ be the Spin$^c$ structure associated to an
intersection point $\mathbf y$, $\widehat{F'}$ be the surface in
$Y_0$ obtained by capping off the boundary of $F'$.

We want to compute $\langle c_1(\underline{\mathfrak
s}_{w,z}(\mathbf y)),[\widehat{F'}]\rangle$.

$(\Sigma,\mbox{\boldmath${\gamma}$},(\mbox{\boldmath${\beta}$}\backslash\{\mu\})\cup\{\lambda''\},w'')$
is a Heegaard diagram for $Y_0$, we wind $\lambda''$ once along
$\delta'\cup\overline{\delta'}$ to create two new intersection
points with $\gamma_1'$. The variant of $\lambda''$ after winding
is denoted by $\lambda^*$. Let $\mathbf y^*$ be an intersection
point close to $\mathbf y$ in this new diagram. A standard
computation of $\langle c_1(\mathfrak s(\mathbf
y^*)),[\widehat{F'}]\rangle$ shows that, an intersection point
$\mathbf x$ is a generator of $\widehat{CFK}(Y,K,-g)$, if and only
if $\mathbf x$ is supported outside $\mathrm{int}(A')$. Now if
$\mathbf x$ is supported outside $\mathrm{int}(A')$, then the
$\beta_2',\dots,\beta_{2g+1+r'}'$ components of $\mathbf x$ have
to lie in $B'$. Hence they are also the
$\gamma'_2,\dots,\gamma'_{2g+1+r'}$ components of $\mathbf x$. So
$\mathbf x$ has no component in $\mathrm{int}(A'')$. This finishes
the proof of the claim.

Using the previous claim, one sees that
$$\widehat{CFK}(Y,K,-g)\cong\widehat{CFK}(Y',K',-g)\otimes\widehat{CFK}(Y'',K'',-g)$$
as abelian groups. Suppose $\Phi$ is a holomorphic disk
 for $\widehat{CFK}(Y,K,-g)$, by the previous claim all the corners of $\Phi$
 are supported outside $A'\cup A''$, so
$\Phi\cap(A'\cup A'')$ is a nonnegative relative periodic domain
in $A'\cup A''$, our previous conclusion before the claim shows
that $\Phi$ is supported away from $\lambda'',\lambda'$. Moreover,
$\Phi$ is supported away from $\gamma'_1\cap\mu$, since $\Phi$
should avoid $w,z$, which lie on different sides of $\mu$. By the
same reason, if $\Phi',\Phi''$ are holomorphic disks for
$\widehat{CFK}(Y',K',-g)$ and $\widehat{CFK}(Y'',K'',-g)$,
respectively, then they are supported away from
$\alpha_1'\cap\mu'$ and $\alpha_1''\cap\mu''$, respectively. Hence
$\Phi$ is the disjoint union of two holomorphic disks for
$\widehat{CFK}(Y',K',-g)$ and $\widehat{CFK}(Y'',K'',-g)$,
respectively. Now our desired result is obvious.
\end{proof}

As a corollary, we have
\begin{cor}
Let $K\subset Y$ be a null-homologous knot with a genus $g$
Seifert surface $F$, $Y_m$ be the $m$--fold cyclic branched cover
of $Y$ over $K$, with respect to $F$, and $K_m$ is the image of
$K$ in $Y_m$. Then
$$\widehat{HFK}(Y_m,K_m,[F],g;\mathbb F)\cong\widehat{HFK}(Y,K,[F],g;\mathbb F)^{\otimes m}$$
as linear spaces over any field $\mathbb F$.\hfill\qedsymbol
\end{cor}

Knot Floer homology of knots in cyclic branched covers has been
 studied by Grigsby \cite{Gr}, with emphasis on
2--bridge knots in $S^3$.

Theorem \ref{HorKnot} can be re-stated in the language of
$\widehat{HFS}(M,\gamma)$ as follows.

\begin{thm}\label{HoriHFS}
Suppose $(M,\gamma)$ is a balanced sutured manifold, and $S\subset
M$ is a horizontal surface. Decompose $(M,\gamma)$ along $S$, we
get two balanced sutured manifolds
$(M_1,\gamma_1),(M_2,\gamma_2)$. Then
$$\widehat{HFS}(M,\gamma)\cong\widehat{HFS}(M_1,\gamma_1)\otimes\widehat{HFS}(M_2,\gamma_2)$$
as linear spaces over any field $\mathbb F$.
\end{thm}

\section{Product decomposition}

In this section, we will study sutured manifold decomposition
along product annuli. We are not able to obtain a formula for
non-separating product annuli, but the formula for separating
product annuli is already enough for many applications.

\begin{thm}\label{SVerHFS}
Suppose $(M,\gamma)$ is a balanced sutured manifold,
$R_{\pm}(\gamma)$ are connected. $\mathcal A\subset M$ is a
separating product annulus, and $\mathcal A$ separates $M$ into
two balanced sutured manifolds $(M_1,\gamma_1)$, $(M_2,\gamma_2)$.

Then we have
$$\widehat{HFS}(M,\gamma)\cong\widehat{HFS}(M_1,\gamma_1)\otimes\widehat{HFS}(M_2,\gamma_2)$$
as vector spaces over any field $\mathbb F$.
\end{thm}

In the first two subsections, we will consider the case that
$\gamma$ has only one component, which lies in $M_2$, and
$M_1=R_1\times[0,1]$, where $R_1$ is a compact genus--1 surface
with one boundary component.

\subsection{A Heegaard diagram related to $(M,\gamma)$}

\begin{construction}\label{RelMorse1} Let $\psi\co R_+(\gamma)\to R_-(\gamma)$ be a homeomorphism, such
that $\psi(R_+(\gamma_i))=R_-(\gamma_i)$, $i=1,2$, and
$\psi|R_+(\gamma_1)$ maps $x\times1$ to $x\times0$ for any $x\in
R_1$. If we glue $R_+(\gamma)$ to $R_-(\gamma)$ by $\psi$, then we
get a manifold with boundary consisting of a torus. This manifold
can be viewed as the complement of a knot $K$ in a manifold $Y$.
We will construct a Heegaard diagram for the pair $(Y,K)$. The
construction is similar to Construction~\ref{RelMorse}.

\noindent{\bf Step 0.}\quad{\sl A relative Morse function}

Consider a self-indexed relative Morse function $u$ on
$(M_2,\gamma_2)$. Let $\widetilde F_2=u^{-1}(\frac32)$.
$\widetilde F_2$ has two boundary components. We denote the one
that lies in the separating annulus $\mathcal A$ by $\widetilde
a$. The other boundary component is denoted by
$\widetilde{\lambda}$. Similarly, the boundary components of
$R_{\pm}(\gamma_2)$ are denoted by $a_{\pm},\lambda_{\pm}$.

Suppose the genus of $R_+(\gamma_2)$ is $g-1$, and $u$ has $r$
index--$1$ critical points, then the genus of $\widetilde F_2$ is
$g+r-1$. The gradient $-\nabla u$ generates a flow $\phi_t$ on
$M_2$.

\noindent{\bf Step 1.}\quad{\sl Construct curves for
$(M_2,\gamma_2)$}

Choose
$$\widetilde{D},\widetilde{\delta}\subset\widetilde{F}_2,\quad
D_{\pm},\delta_{\pm}\subset R_{\pm}(\gamma_2)$$ as in Construction
\ref{RelMorse}. Let
$B_{\pm}=R_{\pm}(\gamma_2)-\mathrm{int}(D_{\pm})$, $\widetilde
B=\widetilde{F}_2-\mathrm{int}(\widetilde{D})$.

On $\widetilde F_2$, there are simple closed curves
$\widetilde{\alpha}_{2g+2},\dots,\widetilde{\alpha}_{2g+1+r}$,
$\widetilde{\beta}_{2g+2},\dots,\widetilde{\beta}_{2g+1+r}$, which
correspond to the critical points of $u$.

The $\widetilde{\alpha}$--curves do not separate $\widetilde F_2$,
so there is an arc $\widetilde{\sigma}\subset\widetilde B$
connecting $\widetilde{\lambda}$ to $\widetilde a$, and
$\widetilde{\sigma}$ is disjoint from $\widetilde{\delta}$ and
$\widetilde{\alpha}$--curves. Similarly, there is an arc
$\widetilde{\tau}\subset\widetilde F_2$ connecting
$\partial\widetilde D$ to $\widetilde a$, and $\widetilde{\tau}$
is disjoint from $\widetilde{\delta}$ and
$\widetilde{\beta}$--curves. Moreover, by stabilization as shown
in Figure~\ref{FigStab}, we can assume
$\widetilde{\sigma}\cap\widetilde{\tau}=\emptyset$, and
$\widetilde{\sigma}$ intersects exactly one
$\widetilde{\beta}$--curve transversely once. Suppose this curve
is $\widetilde{\beta}_{2g+2}$. Let $\sigma_-\subset B_-$ be the
image of $\widetilde{\sigma}$ under the flow $\phi_t$, and
$\sigma_+=\psi^{-1}(\sigma_-)$. Let $\tau_+\subset B_+$ be the
image of $\widetilde{\tau}$ under the flow $\phi_{-t}$.

\begin{figure}
\begin{center}
\begin{picture}(375,120)

\put(70,0){\scalebox{0.4}{\includegraphics*[5pt,365pt][585pt,
660pt]{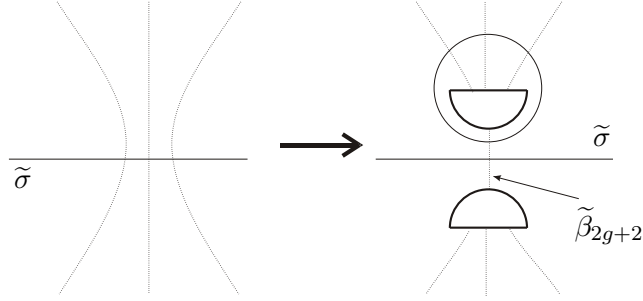}}}

\put(75,45){$\widetilde{\sigma}$}

\put(294,61){$\widetilde{\sigma}$}

\put(287,27){$\widetilde{\beta}_{2g+2}$}

\end{picture}
\caption{A stabilization will eliminate the extra intersection
points of $\widetilde{\sigma}$ with $\widetilde{\tau}$ and
$\widetilde{\beta}$--curves. On the right hand side the two
semicircular holes are glued together.}\label{FigStab}
\end{center}
\end{figure}

Choose $2g-2$ disjoint arcs $\xi^-_4,\dots,\xi^-_{2g+1}\subset
B_-$, such that their endpoints lie on $\lambda_-$, and they are
linearly independent in $\homo_1(B_-,\partial B_-)$. We also
suppose they are disjoint from $\delta_-,\sigma_-$ and the bad
points. Let $\xi^+_i=\psi^{-1}(\xi^-_i)$. We also flow back
$\xi^-_4,\dots,\xi^-_{2g+1}$ by $\phi_{-t}$ to $\widetilde F_2$,
the images are denoted by
$\widetilde{\xi}_4,\dots,\widetilde{\xi}_{2g+1}$.

Choose $2g-2$ disjoint arcs $\eta^+_4,\dots,\eta^+_{2g+1}\subset
B_+$, such that their endpoints lie on $\partial D_+$, and they
are linearly independent in $\homo_1(B_+,\partial B_+)$. We also
suppose they are disjoint from $\delta_+,\tau_+$ and the bad
points. Flow them by $\phi_t$ to $\widetilde F$, the images are
denoted by $\widetilde{\eta}_4,\dots,\widetilde{\eta}_{2g+1}$.

We can slide $\widetilde{\eta}$--curves over
$\widetilde{\beta}_{2g+2}$ to eliminate the possible intersection
points between $\widetilde{\eta}$--curves and
$\widetilde{\sigma}$.

By stabilization, we can assume $\widetilde{\tau}$ does not
intersect $\widetilde{\xi}$--curves, and it intersects exactly one
$\widetilde{\alpha}$--curve transversely once. This curve is
denoted by $\widetilde{\alpha}_{2g+2}$.

\begin{figure}
\begin{center}
\begin{picture}(375,238)

\put(0,0){\scalebox{0.67}{\includegraphics*[5pt,390pt][560pt,
645pt]{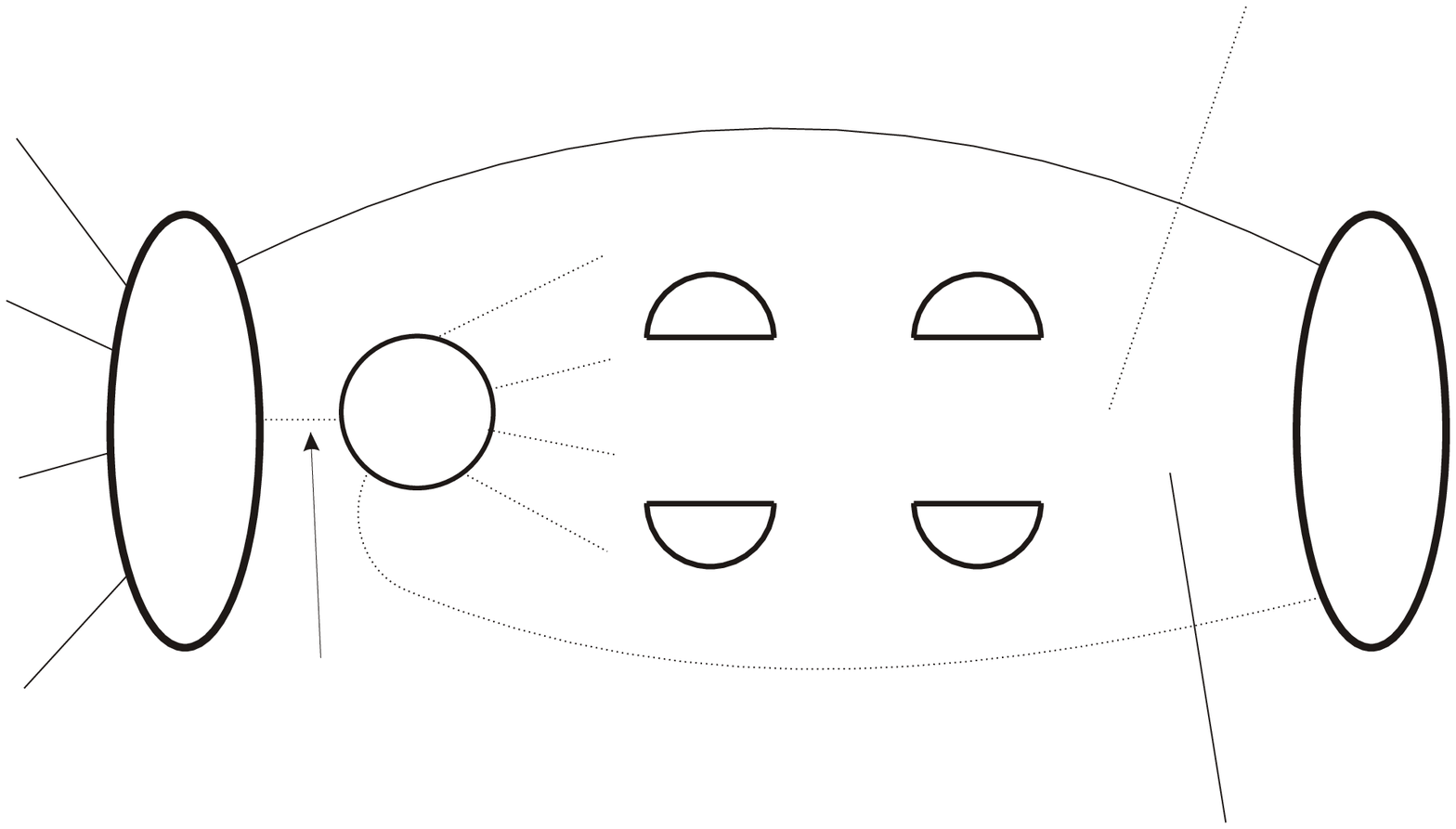}}}

\put(88,148){$\widetilde{\sigma}$}

\put(79,19){$\widetilde{\delta}$}

\put(56,108){$\widetilde{\lambda}$}

\put(104,90){$\widetilde D$}

\put(338,100){$\widetilde a$}

\put(254,21){$\widetilde{\tau}$}

\put(294,107){$\widetilde{\beta}_{2g+2}$}

\put(314,12){$\widetilde{\alpha}_{2g+2}$}

\put(22,152){$\widetilde{\xi}_{4}$}

\put(140,135){$\widetilde{\eta}_4$}

\end{picture}
\caption{A schematic picture of $\widetilde F_2$. The pairs of
semicircular holes are glued together.}
\end{center}
\end{figure}

\noindent{\bf Step 2.}\quad{\sl Find geometric duals of
$\widetilde{\alpha}$--curves}

As in Lemma \ref{admissible}, we can prove that
$$E=\widetilde{B}-\cup_{i=4}^{2g+1}\widetilde{\xi}_i-\cup_{j=2g+2}^{2g+1+r}\widetilde{\alpha}_j-\widetilde{\delta}-\widetilde{\sigma}$$
is connected. Then there are simple closed curves
$$\theta_{2g+2},\dots,\theta_{2g+1+r}\subset\widetilde B-\cup\widetilde{\xi}_i-\widetilde{\sigma}-\widetilde{\delta},$$ which
are geometrically dual to
$\widetilde{\alpha}_{2g+2},\dots,\widetilde{\alpha}_{2g+1+r}$.

We can slide
$\widetilde{\theta}_{2g+3},\dots,\widetilde{\theta}_{2g+1+r}$ over
$\widetilde{\alpha}_{2g+2}$ to eliminate the possible intersection
points between
$\widetilde{\theta}_{2g+3},\dots,\widetilde{\theta}_{2g+1+r}$ and
$\widetilde{\tau}$.

\noindent{\bf Step 3.}\quad{\sl Construct curves for
$(M_1,\gamma_1)$}

As in Figure~\ref{FigPants}, we choose 4 properly embedded arcs
$\xi_2,\xi_3,\eta_2,\eta_3$ on $R_1$. $\xi_2$ intersects $\eta_2$
at one point, and $\xi_3$ intersects $\eta_3$ at one point. There
are no more intersection points between these arcs.

Let $$\widetilde A=\widetilde B\cup_{\widetilde a=(\partial R_1
)\times\frac12}(R_1\times\frac12),$$
$$A_+=B_+\cup_{a_+=(\partial R_1 )\times1}(R_1\times1).$$
Take 4 parallel copies of $\widetilde{\sigma}$, glue them to
$\xi_2\times\frac12,\xi_3\times\frac12$, we get two arcs
$\widetilde{\xi}_2,\widetilde{\xi}_3\subset\widetilde A$. Glue 4
parallel copies of $\widetilde{\tau}$ with
$\eta_2\times\frac12,\eta_3\times\frac12$, we get two arcs
$\widetilde{\eta}_2,\widetilde{\eta}_3\subset\widetilde A$.
Similarly, we can construct
$\xi^+_2,\xi^+_3,\eta^+_2,\eta^+_3\subset A_+$.

\begin{figure}
\begin{center}
\begin{picture}(375,202)
\put(0,0){\scalebox{0.65}{\includegraphics*[1pt,250pt][571pt,
560pt]{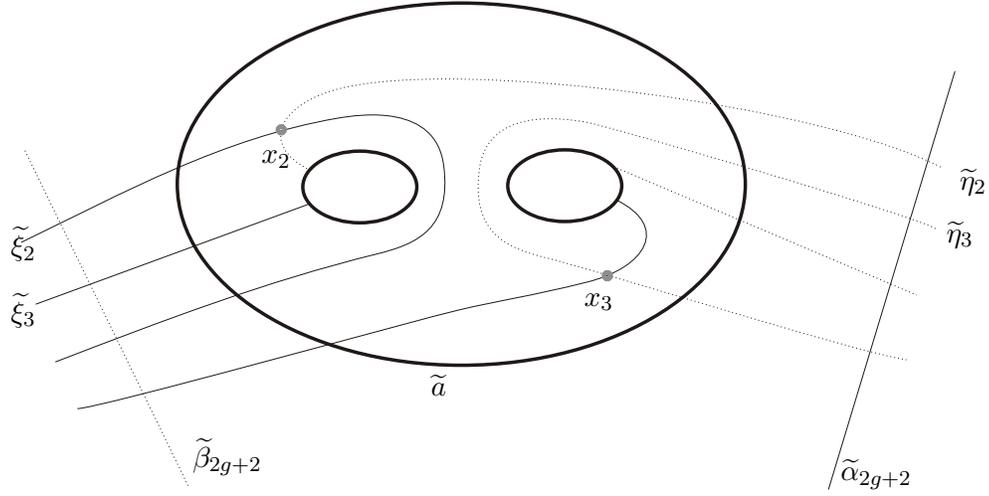}}}

\put(1,100){$\widetilde{\xi}_2$}

\put(1,75){$\widetilde{\xi}_3$}

\put(360,125){$\widetilde{\eta}_2$}

\put(355,105){$\widetilde{\eta}_3$}

\put(96,135){$x_2$}

\put(218,81){$x_3$}

\put(160,47){$\widetilde a$}

\put(70,20){$\widetilde{\beta}_{2g+2}$}

\put(315,15){$\widetilde{\alpha}_{2g+2}$}

\end{picture}
\caption{Local picture of $(\Sigma,\mbox{\boldmath${\alpha}$},
\mbox{\boldmath$\beta$},w,z)$ near $R_1\times\frac12$. The surface
$R_1$ is obtained by gluing the  boundaries of two holes in the
disk by a reflection.}\label{FigPants}
\end{center}
\end{figure}

\noindent{\bf Step 4.}\quad{\sl Construct a Heegaard diagram}

Let
$$\Sigma=A_+\cup \widetilde A\cup\{\textrm{2 annuli}\},$$
where one annulus is $\phi_{[0,3/2]}(\partial D_+)$, the other is
$\phi_{[0,3/2]}(\widetilde{\lambda})$. Construct
$\alpha_i,\beta_i$, ($i=2,\dots,2g+1$) $\alpha_1,\lambda,\mu$ as
in Step 2 of Construction~\ref{RelMorse}.

Let
\begin{eqnarray*}
\mbox{\boldmath${\alpha}$} &= &
\{\alpha_2,\dots,\alpha_{2g+1}\}\cup\{\widetilde{\alpha}_{2g+2},\dots,\widetilde{\alpha}_{2g+1+r}\}\cup\{\alpha_1\},\\
\mbox{\boldmath${\beta}$} &= &
\{\beta_2,\dots,\beta_{2g+1}\}\cup\{\widetilde{\beta}_{2g+2},\dots,\widetilde{\beta}_{2g+1+r}\}\cup\{\mu\}.
\end{eqnarray*}
Pick two base points $w,z$ near $\lambda\cap\mu$, but on different
sides of $\mu$. As in Construction~\ref{RelMorse},
$$(\Sigma,\mbox{\boldmath${\alpha}$},
\mbox{\boldmath$\beta$},w,z)$$ is a Heegaard diagram for $(Y,K)$.
\qed
\end{construction}

It is easy to check that the Heegaard diagram constructed above is
a sutured Heegaard diagram. In order to prove our desired result,
we still need to change the diagram by handleslides.

\begin{figure}
\begin{center}
\begin{picture}(375,205)
\put(0,0){\scalebox{0.65}{\includegraphics*[1pt,250pt][571pt,
565pt]{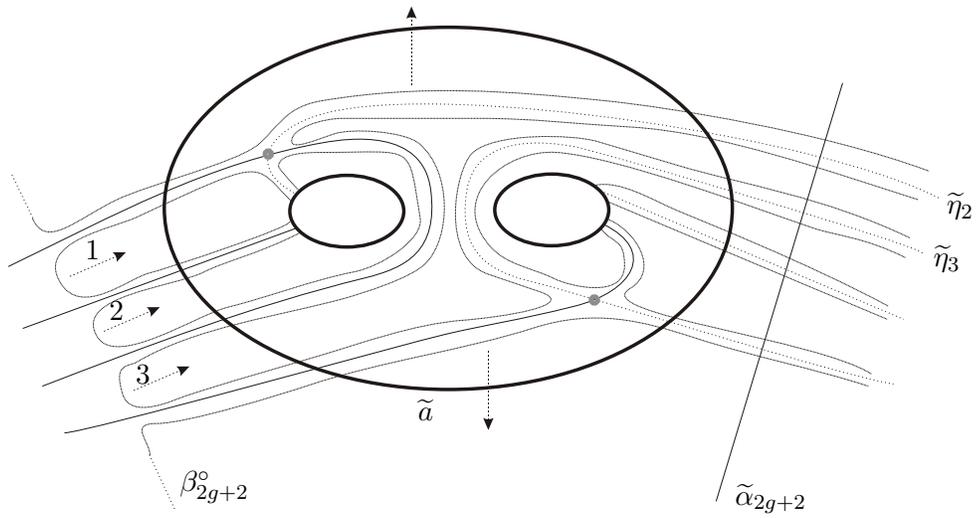}}}

\put(360,125){$\widetilde{\eta}_2$}

\put(355,105){$\widetilde{\eta}_3$}

\put(35,108){1}

\put(44,85){2}

\put(54,61){3}

\put(160,47){$\widetilde a$}

\put(70,20){$\beta^{\circ}_{2g+2}$}

\put(280,15){$\widetilde{\alpha}_{2g+2}$}

\end{picture}
\caption{Local picture of $(\Sigma,\mbox{\boldmath${\alpha}$} ,
\mbox{\boldmath$\beta$}^{\circ},w,z)$ near $R_1\times\frac12$. The
arrows indicate the directions of a further
isotopy.}\label{FigSlide1}
\end{center}
\end{figure}

\begin{figure}
\begin{center}
\begin{picture}(375,205)
\put(0,0){\scalebox{0.65}{\includegraphics*[1pt,250pt][571pt,
565pt]{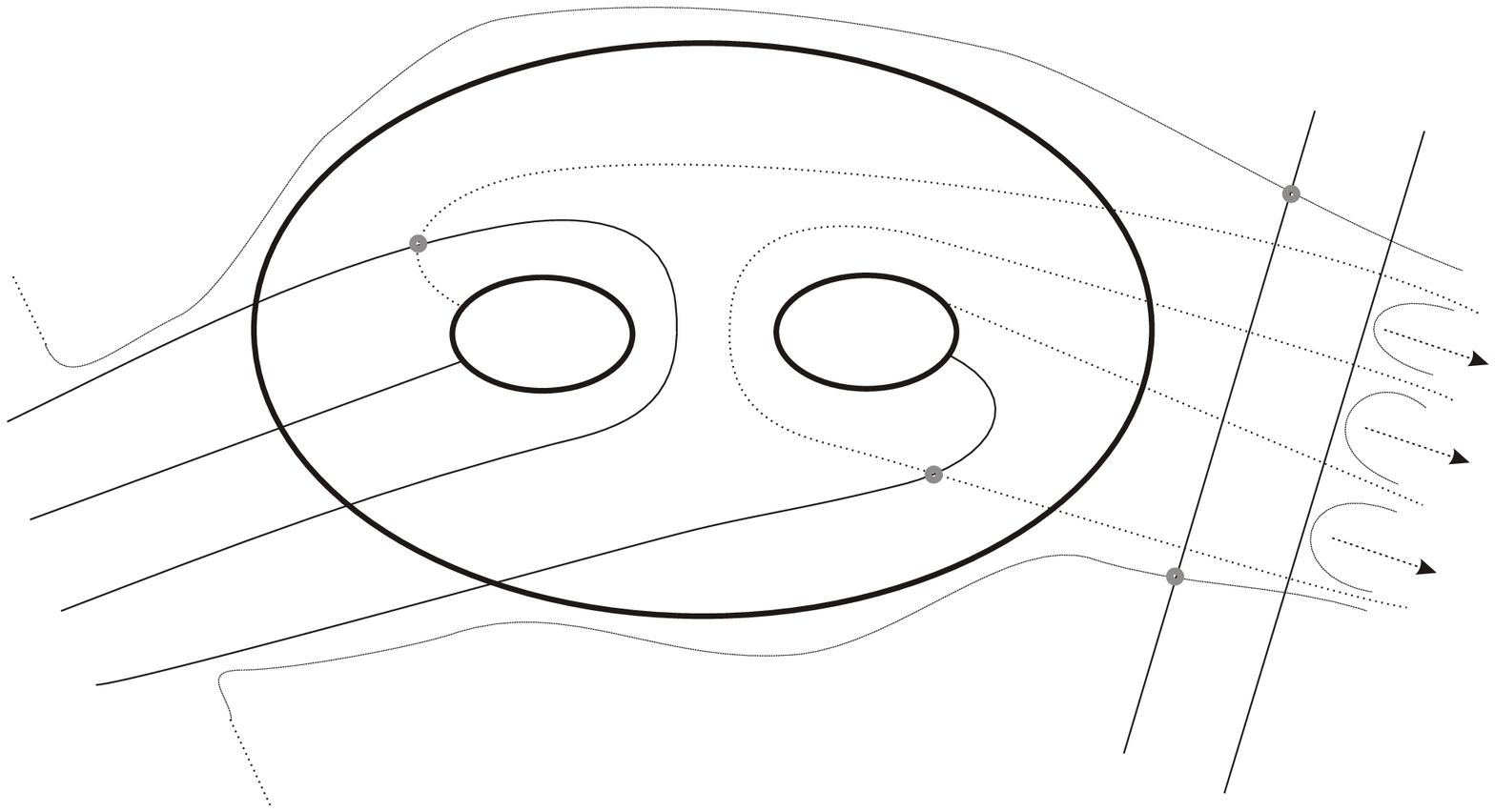}}}

\put(1,100){$\widetilde{\xi}_2$}

\put(1,75){$\widetilde{\xi}_3$}

\put(360,125){$\widetilde{\eta}_2$}

\put(355,105){$\widetilde{\eta}_3$}

\put(96,135){$x_2$}

\put(218,100){$x_3$}

\put(70,20){$\beta^*_{2g+2}$}

\put(255,15){$\widetilde{\alpha}_{2g+2}$}

\put(300,14){$\alpha_1$}

\put(330,75){1}

\put(332,100){2}

\put(339,123){3}

\put(304,170){$x$}

\put(272,61){$y$}
\end{picture}
\caption{Local picture of $(\Sigma,\mbox{\boldmath${\alpha}$} ,
\mbox{\boldmath$\beta$}^*,w,z)$.}\label{FigSlide2}
\end{center}
\end{figure}

\begin{construction}\label{slide}
In the Heegaard diagram $(\Sigma,\mbox{\boldmath${\alpha}$},
\mbox{\boldmath$\beta$},w,z)$ constructed in
Construction~\ref{RelMorse1}, the curve $\widetilde{\beta}_{2g+2}$
has 4 intersection points with $\widetilde{\xi}_2$ and
$\widetilde{\xi}_3$. Since $\widetilde{\xi}_2$ intersects
$\beta_2$ exactly once, we can slide $\widetilde{\beta}_{2g+2}$
over $\beta_2$ twice, to eliminate the intersection points between
$\widetilde{\beta}_{2g+2}$ and $\widetilde{\xi}_2$. Note that in
the two handleslides, the orientations of $\beta_2$ are different,
so the new curve is homologous to $\widetilde{\beta}_{2g+2}$.
Similarly, we can slide $\widetilde{\beta}_{2g+2}$ over $\beta_3$
twice, to eliminate the intersection points between
$\widetilde{\beta}_{2g+2}$ and $\widetilde{\xi}_3$. The local
picture is shown in Figure~\ref{FigSlide1}. The new
$\widetilde{\beta}$--curve is denoted by $\beta^{\circ}_{2g+2}$.

Moreover, we can isotope $\beta^{\circ}_{2g+2}$ to eliminate its
intersection points with $\widetilde{\alpha}_{2g+2}$ and
$\alpha_1$, then isotope it so that it lies outside
$R_1\times\frac12$, as shown in Figure~\ref{FigSlide1} and
Figure~\ref{FigSlide2}. The new $\widetilde{\beta}$--curve is
denoted by $\beta^*_{2g+2}$. $\beta^*_{2g+2}$ is homologous to
$\widetilde{\beta}_{2g+2}$ in $\homo_1(\Sigma)$, thus
$\beta^*_{2g+2}\cap A_+$ is null-homologous in
$\homo_1(A_+,\partial A_+)$.

The new Heegaard diagram after handlesliding is denoted by
$(\Sigma,\mbox{\boldmath${\alpha}$},
\mbox{\boldmath$\beta$}^*,w,z)$. This diagram is still a sutured
Heegaard diagram. \qed
\end{construction}

\begin{lem}\label{support}
After winding transverse to $\widetilde{\alpha}$--curves in
$\widetilde B$, and transverse to $\xi^+$--curves in $A_+$, we can
get a weakly admissible Heegaard diagram for $(Y,K)$. Moreover, in
this diagram, any holomorphic disk for $\widehat{CFK}(Y,K,-g)$ is
supported in $\widetilde B$.
\end{lem}
\begin{proof}
In Step 2 of Construction \ref{RelMorse1}, we find circles
$\theta_{2g+2},\dots,\theta_{2g+1+r}$, which are geometrically
dual to
$\widetilde{\alpha}_{2g+2},\dots,\widetilde{\alpha}_{2g+1+r}$.
Moreover, $\theta_{2g+3},\dots,\theta_{2g+1+r}$ are disjoint with
$\widetilde{\xi}$--curves, $\widetilde{\delta}$,
$\widetilde{\sigma}$ and $\widetilde{\tau}$. As in the proof of
\cite[Proposition 3.3]{Ni2}, we can wind
$\widetilde{\alpha}_{2g+3},\dots,\widetilde{\alpha}_{2g+1+r}$
along $\theta_{2g+3},\dots,\theta_{2g+1+r}$, and wind
$\xi^+$--curves in $A_+$, to get a new Heegaard diagram. In this
diagram, if $\mathcal P$ is a nonnegative periodic domain, then
$\partial\mathcal P$ does not contain
$\widetilde{\alpha}_{2g+3},\dots,\widetilde{\alpha}_{2g+1+r}$ and
$\xi^+$--curves. Obviously, $\partial P$ does not contain $\mu$.
Note that $\alpha_1$ is the only attaching curve in the Heegaard
diagram which intersects $\mu$, and $\partial
P\cdot\mu=0\ne\alpha_1\cdot\mu$, so $\partial P$ does not contain
$\alpha_1$.

$\partial(\mathcal P\cap A_+)-\partial A_+$ is a linear combination
of $\beta^*_{2g+2}\cap A_+$ and $\eta^+_i$'s. Note that
$\eta^+_2,\dots,\eta^+_{2g+1}$ are linearly independent in
$\homo_1(A_+,\partial A_+)$, and $\beta^*_{2g+2}\cap A_+$ is
null-homologous in $\homo_1(A_+,\partial A_+)$, so
$\partial(\mathcal P\cap A_+)-\partial A_+$ is a multiple of
$\beta^*_{2g+2}\cap A_+$.

Hence $\partial\mathcal P$ is a linear combination of
$\widetilde{\alpha}_{2g+2}$ and $\widetilde{\beta}$--curves. Here
$\beta^*_{2g+2}$ is also viewed as a $\widetilde{\beta}$--curve,
since it is homologous to $\widetilde{\beta}_{2g+2}$. In Step 1 of
Construction~\ref{RelMorse1}, $\widetilde{\alpha}_{2g+2}$ is
obtained by stabilization. Thus there is a
$\widetilde{\beta}$--curve, which intersects
$\widetilde{\alpha}_{2g+2}$ exactly once.
$\widetilde{\beta}$--curves are mutually disjoint, so
$\widetilde{\alpha}_{2g+2}$ is not homologous to the linear
combination of $\widetilde{\beta}$--curves. This shows that the
new diagram is weakly admissible. This diagram is still denoted by
$(\Sigma,\mbox{\boldmath${\alpha}$},
\mbox{\boldmath$\beta$}^*,w,z)$.

The generators of $\widehat{CFK}(Y,K,-g)$ are supported outside
$\mathrm{int}(A_+)$.

$\alpha_2$ and $\alpha_3$ intersect only two $\beta$--curves
outside $A_+$. Therefore the two intersection points
$x_2=\widetilde{\xi}_2\cap\widetilde{\eta}_2$ and
$x_3=\widetilde{\xi}_3\cap\widetilde{\eta}_3$ must be chosen.

Suppose $\Phi$ is a holomorphic disk for $\widehat{CFK}(Y,K,-g)$.
$\Phi\cap A_+$ is a relative periodic domain \cite[Definition
3.1]{Ni2} in $A_+$. As before, $\Phi$ does not contain any
$\xi_+$--curve after winding, so $\Phi$ is supported away from
$\lambda$. Now $\partial(\Phi\cap A_+)-\partial A_+$ is a linear
combination of $\beta^*_{2g+2}\cap A_+$ and $\eta^+_i$'s. As in
the second paragraph of this proof, we have that
$$\partial(\Phi\cap A_+)-\partial A_+=m(\beta^*_{2g+2}\cap A_+).$$

$\Phi$ is supported away from $\mu$, $\mu$ intersects $\alpha_1$
at exactly one point, so the contribution of $\alpha_1$ to
$\partial\Phi$ is 0. $\beta^*_{2g+2}\cap A_+$ separates $\lambda$
from $\eta_2^+\cap\alpha_1$. The local multiplicities of $\Phi$ in
a small neighborhood of $x_2,x_3$ are 0, and $\Phi$ does not
intersect $\lambda$. As in Figure~\ref{FigDisk}, we conclude that
the local multiplicities of $\Phi$ in two corners at $x$ are 0,
and the local multiplicities of $\Phi$ in two corners at $y$ are
0.

\begin{figure}
\begin{center}
\begin{picture}(320,224)
\put(0,0){\scalebox{0.80}{\includegraphics*[100pt,320pt][500pt,
600pt]{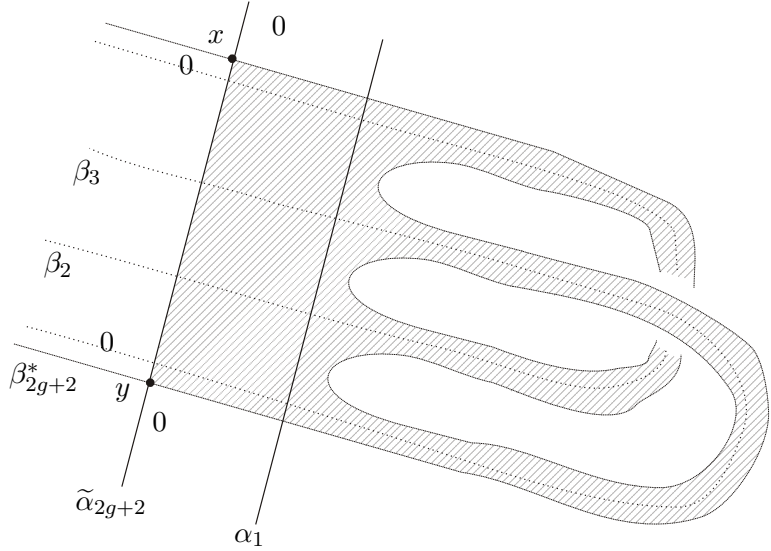}}}

\put(40,20){$\widetilde{\alpha}_{2g+2}$}

\put(100,8){$\alpha_1$}

\put(55,63){$y$}

\put(50,80){0}

\put(70,50){0}

\put(90,197){$x$}

\put(80,185){0}

\put(115,200){0}

\put(15,68){$\beta^*_{2g+2}$}

\put(28,110){$\beta_2$}

\put(39,145){$\beta_3$}
\end{picture}
\caption{A component of the possible holomorphic disk. The two
bands lie in $R_1\times1$.}\label{FigDisk}
\end{center}
\end{figure}

If $m>0$, then $\Phi$ (which is possibly reducible) has a
component $\Phi_0$, which has two vertices $x,y$ and two edges
lying in $\widetilde{\alpha}_{2g+2}$ and $\beta^*_{2g+2}$. See
Figure~\ref{FigDisk} for the picture of $\Phi_0$. The genus of
$\Phi_0$ is 1, and its Maslov index is $-1$. So $\Phi$ has no
contribution to the boundary map in $\widehat{CFK}(Y,K,-g)$.

Hence if $\Phi$ is a holomorphic disk for $\widehat{CFK}(Y,K,-g)$,
then $m=0$, so $\Phi$ is supported in $\widetilde B$.
\end{proof}

\subsection{A Heegaard diagram related to $(M_2,\gamma_2)$}

We also need to construct a Heegaard diagram to compute
$\widehat{HFS}(M_2,\gamma_2)$. As in Subsection~\ref{SecHFS}, we
can add a product 1--handle to $M_2$ with feet at different
components of $\gamma_2$ to get a new sutured manifold
$(M_3,\gamma_3)$. Then we glue $R_+(\gamma_3)$ with
$R_-(\gamma_3)$ to get the complement of a knot $K_3$ in a
manifold $Y_3$. Then
$$\widehat{HFS}(M_2,\gamma_2)=\widehat{HFK}(Y_3,K_3,g).$$ Our next
task is to construct a Heegaard diagram for $(Y_3,K_3)$.

\begin{figure}
\begin{center}
\begin{picture}(375,212)
\put(0,0){\scalebox{0.64}{\includegraphics*[1pt,280pt][581pt,
610pt]{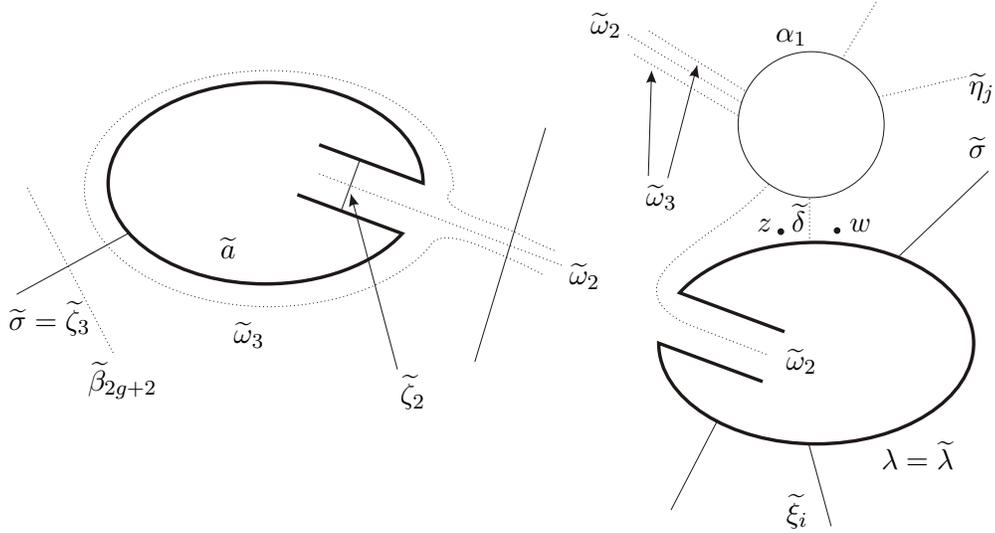}}}

\put(30,60){$\widetilde{\beta}_{2g+2}$}

\put(0,83){$\widetilde{\sigma}=\widetilde{\zeta}_3$}

\put(148,55){$\widetilde{\zeta}_2$}

\put(85,78){$\widetilde{\omega}_3$}

\put(80,110){$\widetilde a$}

\put(212,100){$\widetilde{\omega}_2$}

\put(241,130){$\widetilde{\omega}_3$}

\put(294,68){$\widetilde{\omega}_2$}

\put(283,120){$z$}

\put(296,120){$\widetilde{\delta}$}

\put(318,120){$w$}

\put(290,192){$\alpha_1$}

\put(330,30){$\lambda=\widetilde{\lambda}$}

\put(220,195){$\widetilde{\omega}_2$}

\put(363,148){$\widetilde{\sigma}$}

\put(363,172){$\widetilde{\eta}_j$}

\put(294,10){$\widetilde{\xi}_i$}
\end{picture}
\caption{Local picture of $(\Sigma',\mbox{\boldmath${\alpha}$}',
\mbox{\boldmath$\beta$}',w,z)$, near $\widetilde a$ and near
$\widetilde{\lambda}$. There is a band connecting $\widetilde a$
to $\widetilde{\lambda}$.}\label{FigAdd1}
\end{center}
\end{figure}

\begin{construction}\label{AddBand}
Notations as in Construction \ref{RelMorse1}. As in
Figure~\ref{FigAdd1}, let $\widetilde{C}$ be the surface obtained
by gluing a rectangle $I\times[0,1]$ to $\widetilde B$, such that
$$(I\times[0,1])\cap\widetilde B=I\times\{0,1\},(I\times0)\subset\widetilde{\lambda},(I\times1)\subset\widetilde{a}.$$
Moreover, $\widetilde{\delta}\cap\widetilde{\lambda}$ is not
contained in $I\times0$, but it is close to $I\times0$.
$\widetilde{\tau}\cap\widetilde a\in\mathrm{int}(I\times1)$, and
$\widetilde{\sigma}\cap\widetilde a\notin I\times1$. Let
$\widetilde{\zeta}_2=I\times\frac12$,
$\widetilde{\zeta}_3=\widetilde{\sigma}$. Let
$\widetilde{\omega}_2$ be the union of $\widetilde{\tau}$,
$p\times[0,1]$ ($p\in\mathrm{int}(I)$) and a parallel copy of
$\widetilde{\delta}$. Let $T$ be the neighborhood of
$\widetilde{\tau}\cup\widetilde a$ in $\widetilde B$.
$\widetilde{\omega}_3=T\cap(\widetilde B-\mathrm{int}(T))$.

Similarly, construct a surface $C_+$, and curves
$\zeta^+_2,\zeta^+_3,\omega^+_2,\omega^+_3$ on it.

As in Step 4 of Construction \ref{RelMorse1}, for $i=2,3$, let
$$\alpha'_i=\zeta^+_i\cup\widetilde{\zeta}_i\cup\{\textrm{2
arcs}\},$$
$$\beta'_i=\omega^+_i\cup\widetilde{\omega}_i\cup\{\textrm{2
arcs}\}.$$

Let
\begin{eqnarray*}
\Sigma'&=&C_+\cup \widetilde C\cup\{\textrm{2 annuli}\},\\
\mbox{\boldmath${\alpha}$}' &= &
\{\alpha'_2,\alpha'_3,\alpha_4,\dots,\alpha_{2g+1}\}\cup\{\widetilde{\alpha}_{2g+2},\dots,\widetilde{\alpha}_{2g+1+r}\}\cup\{\alpha_1\},\\
\mbox{\boldmath${\beta}$}' &= &
\{\beta'_2,\beta'_3,\beta_4,\dots,\beta_{2g+1}\}\cup\{\widetilde{\beta}_{2g+2},\dots,\widetilde{\beta}_{2g+1+r}\}\cup\{\mu\}.
\end{eqnarray*}

Then
$$(\Sigma',\mbox{\boldmath${\alpha}$}',
\mbox{\boldmath$\beta$}',w,z)$$ is a Heegaard diagram for
$(Y_3,K_3)$.\qed
\end{construction}

\begin{figure}
\begin{center}
\begin{picture}(375,213)
\put(0,0){\scalebox{0.71}{\includegraphics*[1pt,280pt][526pt,
580pt]{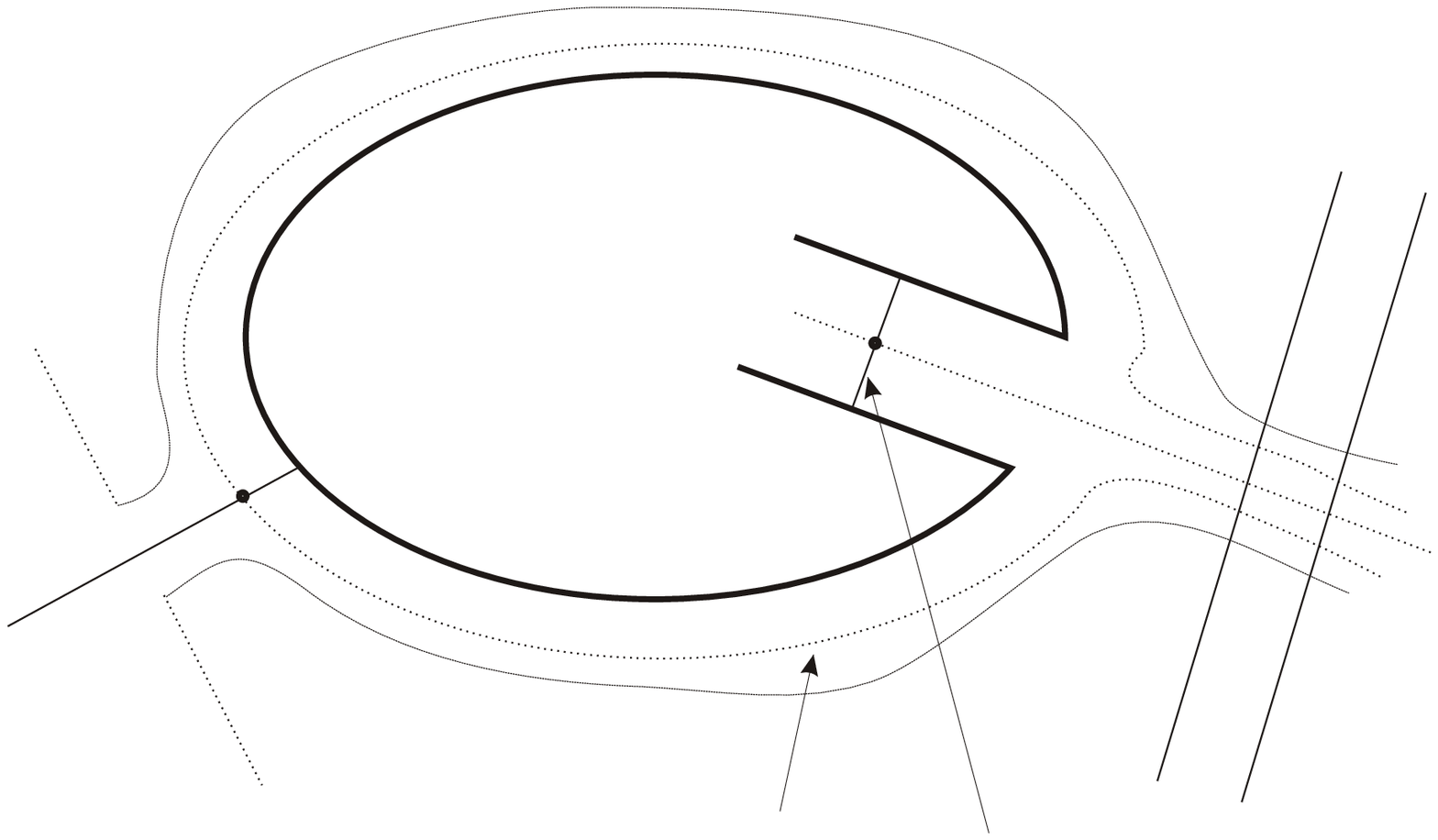}}}

\put(60,9){$\beta''_{2g+2}$}

\put(1,45){$\widetilde{\sigma}=\widetilde{\zeta}_3$}

\put(51,81){$y_3$}

\put(221,131){$y_2$}

\put(246,8){$\widetilde{\zeta}_2$}

\put(195,7){$\widetilde{\omega}_3$}

\put(357,75){$\widetilde{\omega}_2$}
\end{picture}
\caption{Local picture of $(\Sigma',\mbox{\boldmath${\alpha}$}' ,
\mbox{\boldmath$\beta$}'',w,z)$ near $\widetilde
a$.}\label{FigAdd2}
\end{center}
\end{figure}

\begin{construction}
In the diagram $(\Sigma',\mbox{\boldmath${\alpha}$}',
\mbox{\boldmath$\beta$}',w,z)$, we can slide
$\widetilde{\beta}_{2g+2}$ over $\beta'_3$ once, to get a new
curve $\beta''_{2g+2}$. The new diagram is denoted by
$(\Sigma',\mbox{\boldmath${\alpha}$}',
\mbox{\boldmath$\beta$}'',w,z)$. See Figure~\ref{FigAdd2} for the
local picture.\qed
\end{construction}

\begin{lem}\label{support1}
After winding transverse to $\widetilde{\alpha}$--curves in
$\widetilde B$, and transverse to $\xi^+$--curves in $B_+$, we can
get a weakly admissible Heegaard diagram for $(Y_3,K_3)$. In this
diagram, any holomorphic disk for $\widehat{CFK}(Y_3,K_3,-g)$ is
supported in $\widetilde B$.
\end{lem}
\begin{proof}
The proof of this lemma is similar to the proof of
Lemma~\ref{support}, and the argument here is even simpler. Again,
we can wind some curves, such that the diagram becomes weakly
admissible, and the holomorphic disks for
$\widehat{CFK}(Y_3,K_3,-g)$ are supported away from $\lambda'$.
Here $\lambda'\subset\partial\widetilde C$ is the connected sum of
$\lambda$ and $\widetilde a$. Suppose $\Phi'$ is a holomorphic
disk for $\widehat{CFK}(Y_3,K_3,-g)$, then $\Phi'\cap C_+$ is a
relative periodic domain in $C_+$. Since $\Phi'$ is disjoint with
$\lambda'$, $\partial\Phi'\cap C_+$ does not contain
$\xi^+$--curves, and the local multiplicities of $\Phi'$ near
$y_2$ and $y_3$ are zero. So $\partial(\Phi'\cap C_+)-\partial
C_+$ is a linear combination of
$\eta^+_4,\dots,\eta^+_{2g+1},\beta''_{2g+2}\cap C_+$. These arcs
are linearly independent in $\homo_1(C_+,\partial C_+)$, so
$\Phi'$ is supported in $\widetilde B$.
\end{proof}

\subsection{Proof of the product decomposition formula}

\begin{proof}[Proof of Theorem \ref{SVerHFS}]
The surfaces
$$F_{\pm}=R_{\pm}(\gamma_1)\cup\mathcal A\cup
R_{\mp}(\gamma_2)$$ are homologous to $R_+(\gamma)$. $F_{\pm}$ can
be isotoped to horizontal surfaces in $(M,\gamma)$.

$F_+$ decomposes $(M,\gamma)$ into two sutured manifolds
$(M'_1,\gamma'_1)$, $(M'_2,\gamma'_2)$. Here $(M'_1,\gamma'_1)$ is
the sutured manifold obtained by gluing $R_-(\gamma_2)\times I$ to
$(M_1,\gamma_1)$ along $\mathcal A$, and $(M'_2,\gamma'_2)$ is the
sutured manifold obtained by gluing $R_+(\gamma_1)\times I$ to
$(M_2,\gamma_2)$ along $\mathcal A$. Now one can apply Theorem
\ref{HoriHFS} to conclude that
$$\widehat{HFS}(M,\gamma)\cong\widehat{HFS}(M'_1,\gamma'_1)\otimes\widehat{HFS}(M'_2,\gamma'_2).$$
We only need to show that
$\widehat{HFS}(M'_i,\gamma'_i)\cong\widehat{HFS}(M_i,\gamma_i)$
for $i=1,2$. Hence we can reduce our theorem to the case that one
of the two sutured submanifolds $M_1,M_2$ is a product.

From now on, we assume $M_1$ is a product.

According to Definition \ref{DefBal}, $\gamma\ne\emptyset$. If
$M_1\cap\gamma\ne\emptyset$, then one can decompose $M$ along
product disks to get $M_2$. Now we apply
Proposition~{\ref{WellHFS} (2)} to conclude that
$\widehat{HFS}(M,\gamma)\cong\widehat{HFS}(M_2,\gamma_2)$.

Now we consider the case that $M_1\cap\gamma=\emptyset$. By adding
product 1--handles with feet at $\gamma$, we can get a sutured
manifold with connected suture.

$R_+(\gamma_1)$ contains a subsurface $G$ which is a
once-punctured torus. $\partial G\times I$ splits $M$ into two
sutured manifolds $G\times I$ and $(M^*,\gamma^*)$. One can then
decompose $M^*$ along product disks to get $(M_2,\gamma_2)$. Hence
we only need to prove the decomposition formula for the case of
splitting along $\partial G\times I$. From now on, we focus on
this case, namely, the case that the genus of $R_+(\gamma_1)$ is
1.

We apply the constructions in the previous two subsections to get
Heegaard diagrams for $(Y,K)$ and $(Y_3,K_3)$. See
Figure~\ref{FigSlide2} and Figure~\ref{FigAdd2} for the local
pictures. For generators of $\widehat{CFK}(Y,K,-g)$, the two
intersection points $x_2,x_3$ in Figure~\ref{FigSlide2} must be
chosen; for generators of $\widehat{CFK}(Y_3,K_3,-g)$, the two
intersection points $y_2,y_3$ in Figure~\ref{FigAdd2} must be
chosen. Thus the generators of $\widehat{CFK}(Y,K,-g)$ and
$\widehat{CFK}(Y_3,K_3,-g)$ are in one-to-one correspondence. By
Lemma~\ref{support} and Lemma~\ref{support1}, the holomorphic
disks for these two chain complexes are also the same. Hence
$$\widehat{HFK}(Y,K,-g)=\widehat{HFK}(Y_3,K_3,-g),$$ which means
that
$$\widehat{HFS}(M,\gamma)=\widehat{HFS}(M_2,\gamma_2)$$
by definition.
\end{proof}

\subsection{An application to satellite knots}

As an application, we can compute the topmost terms in the knot
Floer homology of satellite knots with nonzero winding numbers. We
recall the following definition from \cite{Ni2}.

\begin{defn}
Suppose $K$ is a null-homologous knot in $Y$, $F$ is a Seifert
surface of $K$ (not necessarily of minimal genus). $V$ is a
3--manifold, $\partial V=T^2$, $L\subset V$ is a nontrivial knot.
$G\subset V$ is a compact connected oriented surface so that $L$
is a component of $\partial G$, and $\partial G-L$ (may be empty)
consists of parallel essential circles on $\partial V$.
Orientations on these circles are induced from the orientation on
$G$, we require that these circles are parallel as oriented ones.
We glue $V$ to $Y-\mathrm{int}(\mathrm{Nd}(K))$, so that any
component of $\partial G-L$ is null-homologous in
$Y-\mathrm{int}(\mathrm{Nd}(K))$. The new manifold is denoted by
$Y^*$, and the image of $L$ in $Y^*$ is denoted by $K^*$. We then
say $K^*$ is a {\it satellite knot} of $K$, and $K$ a {\it
companion knot} of $K^*$. Let $p$ denote the number of components
of $\partial G-L$, $p$ will be called the {\it winding number} of
$L$ in $V$.
\end{defn}

Suppose $p>0$, $F$ is a minimal genus Seifert surface for $K$,
then a minimal genus Seifert surface $F^*$ for $K^*$ can be
obtained as follows: take $p$ parallel copies of $F$, and glue
them to a certain surface $G$ in $V-\mathrm{int}(\mathrm{Nd}(L))$.
We decompose $V-\mathrm{int}(\mathrm{Nd}(L))$ along $G$, the
resulting sutured manifold is denoted by $(M(L),\gamma(L))$, where
$\gamma(L)$ consists of $p+1$ annuli, $p$ of them lie on $\partial
V$, denoted by $A_1,\dots,A_p$.

\begin{cor}
With notations as above, suppose the genus of $K$ is $g>0$, and
the genus of $K^*$ is $g^*$, then
$$\widehat{HFK}(Y^*,K^*,[F^*],g^*)\cong\widehat{HFK}(Y,K,[F],g)\otimes\widehat{HFS}(M(L),\gamma(L))$$
as linear spaces over any field $\mathbb F$.
\end{cor}
\begin{proof}
Let $(M,\gamma)$ be the sutured manifold obtained by decomposing
$Y^*-\mathrm{int}(\mathrm{Nd}(K^*))$ along $F^*$. Note that
$A_1,\dots,A_p$ are separating product annuli in $M$. The desired
result holds by Theorem~\ref{SVerHFS}.
\end{proof}

Matthew Hedden also got some interesting results regarding knot
Floer homology of satellite knots with nonzero winding numbers
\cite{H}. Our result can be compared with his.

\section{Characteristic product regions}

\begin{defn}
Suppose $(M,\gamma)$ is an irreducible sutured manifold, $\gamma$
has no toral component, $R_-(\gamma)$, $R_+(\gamma)$ are
incompressible and diffeomorphic to each other. A {\it product
region} for $M$ is a submanifold $\Phi\times I$ of $N$, such that
$\Phi$ is a compact, possibly disconnected, surface, and
$\Phi\times0$, $\Phi\times1$ are incompressible subsurfaces of
$R_-(\gamma)$, $R_+(\gamma)$, respectively.

There exists a product region $E\times I$, such that if
$\Phi\times I$ is any product region for $M$, then there is an
ambient isotopy of $M$ which takes $\Phi\times I$ into $E\times
I$. $E\times I$ is called a {\it characteristic product region}
for $M$.
\end{defn}

The theory of characteristic product regions is actually a part of
JSJ theory (\cite{JS},\cite{Jo}), the version that we need in the
current paper can be found in \cite{CL}.

The following theorem can be abbreviated as: if
$\widehat{HFS}(M,\gamma)\cong\mathbb Z$, then the characteristic
product region carries all the homology. A version of this theorem
is also proved by Ian Agol via a different approach.

\begin{thm}\label{criterion}
Suppose $(M,\gamma)$ is an irreducible balanced sutured manifold,
$\gamma$ has only one component, and $(M,\gamma)$ is vertically
prime. Let $E\times I\subset M$ be the characteristic product
region for $M$.

If $\widehat{HFS}(M,\gamma)\cong\mathbb Z$, then the map
$$i_*\co\homo_1(E\times I)\to\homo_1(M)$$
is surjective.
\end{thm}

By Proposition \ref{HomoProd}, $(M,\gamma)$ is a homology product.

Let $G$ be a genus--1 compact surface with one boundary component.
Glue the two sutured manifolds $(M,\gamma)$ and $G\times I$
together along their vertical boundaries, we get a sutured
manifold $N$ with empty suture. $N$ has two boundary components
$\Sigma=\Sigma_-=R_-(\gamma)\cup (G\times 0)$,
$\Sigma_+=R_+(\gamma)\cup (G\times 1)$. $N$ is also a homology
product, thus there is a natural isomorphism
$$\partial_*\co\homo_2(N,\partial N)\to\homo_1(\Sigma).$$

\begin{rem}\label{foliation}
Since $(M,\gamma)$ is a homology product, we can glue
$R_-(\gamma)$ to $R_+(\gamma)$, so that the resulting manifold is
the complement of a knot $K$ in a homology 3--sphere $Y\not\approx
S^3$. Suppose $J$ is a genus--$1$ fibred knot in a homology sphere
$Z\not\approx S^3$, $G'$ is a fibre. Consider the knot
$K\#J\subset Y\#Z$. Let $Y'_0$ be the manifold obtained by
0--surgery on $K\#J$. Let $F=R_-(\gamma)\subset Y$, $H$ be the
boundary connected sum of $F$ and $G'$, and $\widehat H$ be the
extension of $H$ in $Y'_0$.

If we cut $Y'_0$ open along $\widehat H$, then we get the manifold
$N$. \cite[Theorem 8.9]{G3} shows that $Y'_0$ admits a taut
foliation, with $\widehat H$ as a compact leaf. Thus this
foliation induces a foliation of $N$.
\end{rem}

Assume that the map
$$i_*\co\homo_1(E\times I)\to\homo_1(M)$$
is not surjective. We can find a simple closed curve $\omega\subset
R_-(\gamma)$, such that $[\omega]$ is not in $i_*(\homo_1(E\times
I))$.

Let $\omega_-=\omega\subset\Sigma_-$, and let
$\omega_+\subset\Sigma_+$ be a circle homologous to $\omega$. We
fix an arc $\delta$ connecting $\Sigma_-$ to $\Sigma_+$. Let
$\mathcal S_m(+\omega)$ be the set of properly embedded surfaces
$S\subset N$, such that $\partial S=\omega_-\sqcup(-\omega_+)$,
and the algebraic intersection number of $S$ with $\delta$ is $m$.
Here $-\omega_+$ denotes the curve $\omega_+$, but with opposite
orientation. Similarly, let $\mathcal S_m(-\omega)$ be the set of
properly embedded surfaces $S\subset N$, such that $\partial
S=(-\omega_-)\sqcup\omega_+$, and the algebraic intersection
number of $S$ with $\delta$ is $m$. Let $x(\mathcal
S_m(\pm\omega))$ be the minimal value of $x(S)$ for all
$S\in\mathcal S_m(\pm\omega)$. It is obvious that $$x(\mathcal
S_{m+1}(\pm\omega))\leq x(\mathcal S_{m}(\pm\omega))+ x(\Sigma).$$

The next fact is implicitly contained in \cite[Theorem 3.13]{G1}.

\begin{lem}\label{ExistTaut}
When $m$ is sufficiently large, there exist connected surfaces
$S_1\in\mathcal S_m(+\omega)$ and $S_2\in\mathcal S_m(-\omega)$,
such that they give taut decompositions of $N$.
\end{lem}
\begin{proof}
Let $D(N)$ be the double of $N$ along $\partial N$.
$a=\partial_*^{-1}([\omega])\in\homo_2(N,\partial N)$ is the
homology class whose intersection with $\Sigma$ is $[\omega]$,
$D(a)$ is its double in $\homo_2(D(N))$. There exists $C\geq0$,
such that if $k>C$, then
$x(D(a)+(k+1)[\Sigma])=x(D(a)+k[\Sigma])+x(\Sigma)$. As in the
proof of \cite[Theorem~3.13]{G1}, if $Q$ is a Thurston norm
minimizing surface in the homology class $D(a)+k[\Sigma]$, and
$Q\cap N$ has no disk or sphere components, then $Q\cap N$ gives a
taut decomposition of $N$.

We can do oriented cut-and-paste of $Q$ with copies of $\Sigma$,
to get a new surface $Q'$, such that $Q'\cap N$ has positive
intersection number with $\delta$. Of course, $Q'\cap N$ still
gives a taut decomposition of $N$. The not-so-good thing is that
$\partial (Q'\cap N)$ is not necessarily
$\omega_-\sqcup(-\omega_+)$. What we can do is to apply
\cite[Lemma 0.6]{G2}. Note that in the proof of \cite[Lemma
0.6]{G2}, one gets a new decomposition surface with prescribed
boundary by gluing subsurfaces $W_i$ of $\Sigma_{\pm}$ to the
original decomposition surface. And by \cite[Lemma~3.10]{G1},
$W_i$ has the same orientation as $\Sigma_{\pm}$. So the algebraic
intersection number of this new decomposition surface with
$\delta$ is no less than $(Q'\cap N)\cdot \delta>0$.

Denote the new decomposition surface by $S_0$, $\partial
S_0=\omega_-\sqcup(-\omega_+)$. Suppose $S_1$ is the component of
$S_0$ which contains $\omega_-$. For homological reason, $S_1$
should also contain $-\omega_+$. Thus other components of $S_0$
are closed surfaces which do not separate $\Sigma_-$ from
$\Sigma_+$. Hence the algebraic intersection number of other
components with $\delta$ is 0. $S_1$ also gives a taut
decomposition of $N$, by \cite[Lemma~0.4]{G2}. So $S_1$ is the
surface we need. Similarly, we can prove the result for $\mathcal
S_m(-\omega)$.
\end{proof}

We also need the following key lemma.

\begin{lem}\label{alpha+>0}
For any positive integers $p,q$,
$$x(\mathcal S_p(+\omega))+x(\mathcal
S_q(-\omega))>(p+q)x(\Sigma).$$
\end{lem}

Suppose $S_1\in\mathcal S_p(+\omega),S_2\in\mathcal S_q(-\omega)$.
Isotope $S_1,S_2$ so that they are transverse. Since $N$ is
irreducible and $S_1,S_2$ are incompressible, we can assume
$S_1\cup S_2-S_1\cap S_2$ has no disk components. Perform oriented
cut-and-paste to $S_1,S_2$, we get a closed surface $P\subset
\mathrm{int}(N)$, with $x(P)=x(S_1)+x(S_2)$. $P$ has no sphere
components, otherwise $S_1\cup S_2-S_1\cap S_2$ would have disk
components.

\begin{figure}
\begin{center}
\begin{picture}(375,105)

\put(0,0){\scalebox{0.83}{\includegraphics*[0pt,450pt][450pt,
575pt]{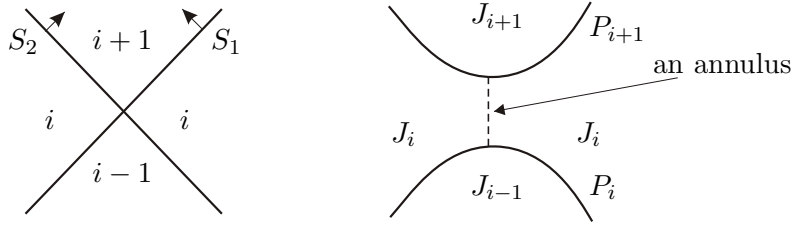}}}

\put(142,75){$S_1$}

\put(65,75){$S_2$}

\put(130,45){$i$}

\put(79,45){$i$}

\put(97,25){$i-1$}

\put(97,75){$i+1$}

\put(310,66){an annulus}

\put(210,40){$J_i$}

\put(280,40){$J_i$}

\put(240,20){$J_{i-1}$}

\put(240,85){$J_{i+1}$}

\put(285,80){$P_{i+1}$}

\put(285,20){$P_i$}

\end{picture}
\caption{Doing oriented cut-and-paste to $S_1$,
$S_2$.}\label{FigCutPaste}
\end{center}
\end{figure}

Now we will deal with the possible toral components of $P$. To
this end, we need the following lemma.

\begin{lem}\label{TorusZero}
If $T\subset\mathrm{int}(N)$ is a torus, then the algebraic
intersection number of $T$ and $\delta$ is 0.
\end{lem}
\begin{proof}
Since $N$ is a homology product, we have
$$\homo_2(D(N))\cong\homo_2(\Sigma)\oplus\homo_1(\Sigma).$$

$T$ is disjoint from $\Sigma$, so $[T]$ must be a multiple of
$[\Sigma]$ in $\homo_2(D(N))$. By Remark \ref{foliation}, $\Sigma$
is Thurston norm minimizing. Since $x(\Sigma)>0=x(T)$, we must
have $[T]=0$. Hence $T\cdot\delta=0$.
\end{proof}

Suppose $T$ is a toral component of $P$, then $T$ is the union of
$2m$ annuli $A_1,A_2,\dots,A_{2m}$, where $A_{2i-1}\subset S_1$,
$A_{2i}\subset S_2$. Let
\begin{eqnarray*}
S_1'&=&(S_1-\bigcup_{i=1}^m A_{2i-1})\cup \bigcup_{i=1}^m (-A_{2i})\\
S_2'&=&(S_2-\bigcup_{i=1}^m A_{2i})\cup \bigcup_{i=1}^m
(-A_{2i-1}).
\end{eqnarray*}
Here $-A_j$ means $A_j$ with opposite orientation.

A small isotopy will arrange that $|S_1'\cap S_2'|<|S_1\cap S_2|$.
Moreover, $x(S_1')=x(S_1)$, $x(S_2')=x(S_2)$. We want to show that
$S_1'\in\mathcal S_p(\omega)$, $S_2'\in\mathcal S_q(-\omega)$.
Obviously, $\partial S_1'=\partial S_1=\omega_-\sqcup(-\omega_+)$,
$\partial S_2'=\partial S_2=(-\omega_-)\sqcup\omega_+$.
Lemma~\ref{TorusZero} shows that $S_1'\cdot\delta=S_1\cdot\delta$.
Thus $S_1'\in\mathcal S_p(+\omega)$. Similarly, $S_2'\in\mathcal
S_q(-\omega)$. Therefore, we can replace $S_1,S_2$ with
$S_1',S_2'$, then continue our argument.

Now we can assume $P$ has no toral components, and proceed to the
proof of Lemma~\ref{alpha+>0}. Our approach to this lemma was
suggested by David Gabai. In fact, this argument is similar to the
argument in \cite[Lemma~8.22]{G3}.

\begin{proof}[Proof of Lemma \ref{alpha+>0}]
If $x(\mathcal S_p(+\omega))+x(\mathcal
S_q(-\omega))\le(p+q)x(\Sigma)$, then we can get a surface
$P\subset \mathrm{int}(N)$ as above, $x(P)\le(p+q)x(\Sigma)$.
Define a function $\varphi\co (N-P)\to\mathbb Z$ as follows. When
$z\in\Sigma_-$, $\varphi(z)=0$. In general, given $z\in N-P$,
choose a path from $\Sigma_-$ to $z$, $\varphi$ is defined to be
the algebraic intersection number of this path with $P$.

$N$ has the homology type of $\Sigma$, thus any closed curve in
$N$ should have zero algebraic intersection number with any closed
surface. Thus $\varphi$ is well-defined. Moreover, the value of
$\varphi$ on $\Sigma_+$ is $p+q$.

Let $J_i$ be the closure of $\{x\in (N-P)|\:\varphi(x)=i\}$,
$P_i=J_{i-1}\cap J_i$. Thus $P=\sqcup_{i=1}^{p+q}P_i$, and
$\cup_{k=0}^{i-1}J_k$ gives a homology between $\Sigma$ and $P_i$.
Since $x(P)\le(p+q)x(\Sigma)$, and $\Sigma$ is Thurston norm
minimizing in $D(N)$, we must have $x(P_i)=x(\Sigma)$ for each
$i$.

$P_i$ has only one component. Otherwise, suppose $P_i=Q_1\sqcup
Q_2$, then $$x(Q_1),x(Q_2)<x(P_i)=x(\Sigma).$$ As in the proof of
Lemma \ref{TorusZero}, we find that $[Q_1],[Q_2]$ are multiples of
$[\Sigma]$, which gives a contradiction.

We can isotope $P$, so that $P_i\cap M$ is a genus $g$ surface
with one boundary component. In fact, after an isotopy, we can
arrange that $P\cap\gamma_1$ consists of parallel essential curves
in $\gamma_1$. Since $F$ and $G$ are Thurston norm minimizing in
$\homo_2(M,\gamma)$ and $\homo_2(G\times I,\partial G\times I)$,
respectively, we must have
$$x(P_i\cap M)=x(F), x(P_i\cap (G\times I))=x(G)=1.$$ If
an annulus $A$ is a component of $P_i\cap(G\times I)$, then we can
isotope $A$ inside $G\times I$ into $\partial G\times I$, a
further isotopy of $P$ will decrease the number of components of
$P\cap(G\times I)$. So we can assume that $P_i\cap(G\times I)$ has
no annular components. Now the fact that $x(P_i\cap (G\times
I))=1$ indicates that $P_i\cap(G\times I)$ is either a
thrice-punctured sphere or a once-punctured torus. In the former
case there would be an essential sphere in $\widehat{G}\times I$,
where $\widehat G$ is the torus obtained by capping off $\partial
G$, which is impossible. Hence $P_i\cap(G\times I)$ is a
once-puncture torus, and $P_i\cap M$ is diffeomorphic to $F$.

Since $M$ is vertically prime, $P_i\cap M$ is parallel to either
$R_-(\gamma)$ or $R_+(\gamma)$ in $M$. Now the picture of
$P_1,\dots,P_{p+q}$ is clear: there exists a number
$r\in\{0,1,\dots,p+q\}$, such that $P_1,\dots,P_r$ are parallel to
$\Sigma_-$, and $P_{r+1},\dots,P_{p+q}$ are parallel to
$\Sigma_+$. Here we let $P_0=\Sigma_-, P_{p+q+1}=\Sigma_+$.

Since $P$ is gotten by doing cut-and-paste to $S_1,S_2$, we can
isotope $S_1$ so that $S_1\cap J_i$ consists of vertical annuli.
See Figure~\ref{FigCutPaste} for the local picture. We denote
$S_1\cap J_i$ by $C_i\times I$, where $C_i$ is the collection of
some circles in $P_i$. Obviously, $[C_i]$ is homologous to
$[\omega]$.

Consider $J_r$, which is bounded by $P_r,P_{r+1}$. $J_r$ is
homeomorphic to $N$. Since
$$[\omega]\notin\mathrm{im}\;(i_*\co\homo_1(E\times I)\to\homo_1(N)),$$
one component of $C_r\times I$ must lie outside
$\mathrm{im}\:i_*$. Thus this vertical annulus can not be
homotoped into $E\times I$, which contradicts to the definition of
characteristic product region.
\end{proof}

\begin{lem}\label{0surgery}Let $K$ be a knot in a homology 3--sphere $Y$,
$Y_p$ be the manifold obtained by $p$--surgery on $K$. Let $g>1$
be the genus of $K$. Suppose $\widehat{HFK}(Y,K,g;\mathbb
Q)\cong\mathbb Q$, then
$$HF^+(Y_0,[g-1];\mathbb Q)\cong\mathbb
Q.$$
\end{lem}
\begin{proof}
We will use $\mathbb Q$ coefficients in the homologies. As in
 \cite[Corollary 4.5]{OSz4}, when $p$ is sufficiently large, we
have two exact triangles
$$\begin{CD}
@>\delta>>\widehat{HFK}(Y,K,g)@>\sigma>> HF^+(Y_p,[g-1]) @>\psi>>
HF^+(Y)@>\delta>>,
\end{CD}
$$
$$\begin{CD}
@>\delta'>>HF^+(Y_0,[g-1])@>\sigma'>> HF^+(Y_p,[g-1]) @>f>>
HF^+(Y)@>\delta'>>.
\end{CD}
$$
And $f$ has the form $\psi+\iota$, where $\iota$ is a sum of
homogeneous maps which have lower orders than $\psi$.

Since $\widehat{HFK}(Y,K,g;\mathbb Q)\cong\mathbb Q$, either
$\delta$ is surjective or $\sigma$ is injective. Therefore, either
$\psi$ is injective, or $\psi$ is surjective. For simplicity,
denote $HF^+(Y_p,[g-1])$ by $A$, and $HF^+(Y)$ by $B$.

If $\psi$ is injective, $B$ can be written as $\psi(A)\oplus C$
for some subgroup $C$ of $B$. If $b\in B$ is in the form of
$(\psi(a),c)$, then let $\rho(b)=a$. Thus $\rho$ is  a
homomorphism, $\rho\psi=\mathrm{id}$, and $\iota\rho\co B\to B$ is
a homomorphism which strictly decreases degree. Now
$$\mathrm{id}-\iota\rho+(\iota\rho)^2-(\iota\rho)^3+\cdots$$
is a well-defined homomorphism. and
$$\psi=(\mathrm{id}-\iota\rho+(\iota\rho)^2-(\iota\rho)^3+\cdots)(\psi+\iota).$$
Hence $f=\psi+\iota$ is also injective, and $HF^+(Y_0,[g-1])\cong
B/f(A)$.

It is easy to check that $\mathrm{id}+\iota\rho$ induces a
homomorphism from $B/\psi(A)$ to $B/f(A)$, whose inverse is
induced by
$$\mathrm{id}-\iota\rho+(\iota\rho)^2-(\iota\rho)^3+\cdots.$$
Thus $B/f(A)\cong B/\psi(A)$. So
$\mathrm{rank}(HF^+(Y_0,[g-1]))=1$.

A similar argument shows that if $\psi$ is surjective, then $f$ is
also surjective, and $\mathrm{rank}(HF^+(Y_0,[g-1]))=1$.
\end{proof}

\begin{proof}[Proof of Theorem \ref{criterion}] Use the notations in
Remark \ref{foliation}, we have
$$\mathrm{rank}(\widehat{HFK}(Y\#Z,K\#J,g+1))=\mathrm{rank}(\widehat{HFK}(Y,K,g))=1.$$
Thus $\mathrm{rank}(HF^+(Y_0',g))=1$ by Lemma~\ref{0surgery}.

If $i_*$ is not surjective, then the proof of \cite[Theorem
1.4]{Gh}, combined with Lemma~\ref{ExistTaut} and
Lemma~\ref{alpha+>0}, shows that
$$\mathrm{rank}(HF^+(Y_0',g))>1,$$
which gives a contradiction.

More precisely, by Lemma~\ref{ExistTaut} and Lemma~\ref{alpha+>0},
there exist connected surfaces $S_1\in\mathcal S_m(+\omega)$ and
$S_2\in\mathcal S_m(-\omega)$, such that they give taut
decompositions of $N$, and $x(S_1)+x(S_2)>2mx(\Sigma)$. By Gabai's
work \cite[Section 5]{G1}, there exist two taut foliations
$\mathscr F_1$, $\mathscr F_2$ of $N$, such that
\begin{eqnarray*}
\chi(S_1)=e(\mathscr F_1,S_1)&=&e(\mathscr F_1,S_0)+m\chi(\Sigma),\\
\chi(S_2)=e(\mathscr F_2,S_2)&=&e(\mathscr
F_2,-S_0)+m\chi(\Sigma).
\end{eqnarray*}
Here $e(\mathscr F,S)$ is defined in \cite[Definition 3.7]{Gh},
$S_0$ is any surface in $\mathcal S_0(+\omega)$.

Now we can conclude that $e(\mathscr F_1,S_0)\ne e(\mathscr
F_2,S_0)$. Hence \cite[Theorem 3.8]{Gh} can be applied.
\end{proof}

\section{Proof of the main theorem}

\begin{proof}[Proof of Theorem \ref{KnotFibre}]
Suppose $(M,\gamma)$ is the sutured manifold obtained by cutting
open $Y-\mathrm{int}(\mathrm{Nd}(K))$ along $F$, $E\times I$ is
the characteristic product region. We need to show that $M$ is a
product. By Proposition~\ref{HomoProd}, $M$ is a homology product.
Moreover, by Theorem~\ref{HorKnot}, we can assume $M$ is
vertically prime.

If $M$ is not a product, then $M-E\times I$ is nonempty. Thus
there exist some product annuli in $(M,\gamma)$, which split off
$E\times I$ from $M$. Let $(M',\gamma')$ be the remaining sutured
manifold. By Theorem~\ref{criterion}, $R_{\pm}(\gamma')$ are
planar surfaces, and $M'\cap(E\times I)$ consists of separating
product annuli in $M$. Since we assume that $M$ is vertically
prime, $M'$ must be connected. (See the first paragraph in the
proof of Theorem~\ref{SVerHFS}.) Moreover, $M'$ is also vertically
prime, and there are no nontrivial product disks or product annuli
in $M'$. By Theorem~\ref{SVerHFS},
$\widehat{HFS}(M',\gamma')\cong\mathbb Z$.

We add some product 1--handles to $M'$ to get a new sutured
manifold $(M'',\gamma'')$ with $\gamma''$ connected. By
Proposition~\ref{WellHFS},
$\widehat{HFS}(M'',\gamma'')\cong\mathbb Z$. It is easy to see
that $M''$ is also vertically prime. Proposition~\ref{HomoProd}
shows that $M''$ is a homology product.

In the manifold $M''$, the characteristic product region
$E''\times I$ is the union of the product 1--handles and
$\mathrm{Nd}(\gamma')$. Obviously
$i_*\co\homo_1(E'')\to\homo_1(M'')$ is not surjective, which
contradicts to Theorem~\ref{criterion}.
\end{proof}

\begin{proof}[Proof of Corollary \ref{LinkFibre}] Cut
$Y-\mathrm{int}(\mathrm{Nd}(L))$ open along $F$, we get a sutured
manifold $(M,\gamma)$, $\widehat{HFS}(M,\gamma)\cong\mathbb Z$. By
adding product 1--handles with feet at $\gamma$, we can get a new
sutured manifold $(M',\gamma')$, where $\gamma'$ has only one
component. We have $\widehat{HFS}(M',\gamma')\cong\mathbb Z$. By
Theorem~\ref{KnotFibre}, $M'$ is a product, hence $M$ is also a
product. So the desired result holds.
\end{proof}

\end{document}